\documentclass[11pt,amstex,amscd]{amsart}

\newcommand{\ot}{\otimes}
\newcommand{\ra}{\rightarrow}

\newcommand{\BC}{\mathbb{C}}
\newcommand{\BR}{\mathbb{R}}
\newcommand{\CP}{\mathbb{CP}}

\newcommand{\BZ}{\mathbb{Z}}
\newcommand{\BH}{\mathbb{H}}
\newcommand{\dirac}{D\!\!\!\!/\,}

\newcommand{\Id}{\mathbf{Id}}
\newcommand{\Imaginary}{\mathrm{Im}}
\newcommand{\euler}{\mathbf{e}}
\newcommand{\I}{\mathbf{i}}

\newcommand{\Spinc}{\ensuremath{\mathrm{Spin}^c}}

\newcommand{\Pic}{\mbox{Pic}}
\newcommand{\Snabla}{\nabla}
\newcommand{\Spnabla}{{\boldsymbol{\nabla}}}
\newcommand{\End}{\mbox{End}}
\newcommand{\T}{{\ensuremath{\textsc{T}}}}
\newcommand{\C}{{\ensuremath{\mathsf{S}}}}

\usepackage{latexsym}
\usepackage{amsmath}
\usepackage{amscd}
\usepackage[all]{xy} 
\usepackage{amsfonts}
\usepackage{amssymb}
\usepackage{graphicx}
\usepackage{amsbsy}

\newtheorem{thm}{Theorem}
\newtheorem{lem}[thm]{Lemma}
\newtheorem{cor}[thm]{Corollary}
\newtheorem{df}[thm]{Definition}
\newtheorem{Example}[thm]{Example}
\newtheorem{prop}[thm]{Proposition}

\newtheorem{remark}[thm]{Remark}

\newtheorem{conj}[thm]{Conjecture}

\newenvironment{rem}{\begin{remark}\rm}{\end{remark}}

\newenvironment{example}{\begin{Example}\rm}{\end{Example}}

\newcounter{results}\stepcounter{results}
\newenvironment{MainThm}{\bigskip \noindent {\bf Theorem \Alph{results}.}
\it \stepcounter{results}}{\bigskip}

\begin{document}


\title{Seiberg-Witten Invariants, Orbifolds, and Circle Actions}

\author{Scott Jeremy Baldridge}
\address{Department of Mathematics, Michigan State University \newline
\hspace*{.375in}East Lansing, Michigan 48824}
\email{\rm{baldridg@math.msu.edu}}
\date{2001}


\begin{abstract}
The main result of this paper is a formula for calculating the
Seiberg-Witten invariants of 4-manifolds with fixed-point free
circle actions.  This is done by showing under suitable
conditions the existence of a diffeomorphism between the moduli
space of the 4-manifold and the moduli space of the quotient
3-orbifold. Two corollaries include $b_+ {>} 1$ $4$-manifolds
with fixed-point free circle actions  are simple type and a new
proof that $\mathcal{SW}_{Y^3\times S^1} = \mathcal{SW}_{Y^3}$.
An infinite number of $b_+{=}1$ $\;4$-manifolds where the
Seiberg-Witten invariants are still diffeomorphism invariants are
constructed and studied.
\end{abstract}

\maketitle

\addtocounter{section}{0}

\section{Introduction}

The main idea of this work is to systematically study 4-manifolds
that admit an $S^1$-action and classify them using Seiberg-Witten
gauge theory.  When the action on $X^4$ is free, the quotient by
the $S^1$-action is a smooth 3-manifold $Y$ and the manifold with
given circle action is classified by the Euler class $\chi \in
H^2(Y;\BZ)$. When the circle action is not free there will be
non-trivial isotropy groups, which forces the orbit space to be an
orbifold rather than a manifold. The main result of this paper is
a formula for calculating the Seiberg-Witten invariants of any
4-manifold with a fixed point free circle action.  We derive the
formula by proving the existence of a diffeomorphism between the
moduli space of the 4-manifold with the moduli space of the
quotient 3-orbifold.

A given manifold may admit more than one circle action. So while
the 3-manifold and Euler class are fully sufficient to classify a
free circle action, the Seiberg-Witten invariants are stronger in
that they are invariant of the underlying space up to
diffeomorphism regardless of the circle action.

The first theorem we prove puts a restriction on the set of
\Spinc\ structures with nontrivial Seiberg-Witten invariants for
manifolds which admit a fixed point free circle action. (See
sections 2 and 3 for descriptions of \Spinc\ structures and
Seiberg-Witten invariants.)

\begin{MainThm}
Let $\xi$ be a \Spinc\ structure on $b_+{\not=}1$ 4-manifold  $X$
with a fixed point free circle action such that
$SW_X(\xi)\not=0$. Then the  \Spinc\ structure $\xi$ is pulled
back from a \Spinc\ structure on $Y$.
\end{MainThm}

See subsection \ref{sec:sw_pullbacks} for the statement when
$b_+(X)=1$.  This theorem is already enough to imply that $X$ is
SW simple type -- that the expected dimension of the moduli space
for all \Spinc\ structures with nontrivial invariants is zero.

Let $\pi:X\ra Y$ be the projection map from a smooth 4-manifold
with a fixed point free circle action to its quotient orbifold.
The manifold $X$ can be thought of as a orbifold circle bundle
over $Y$. If $\I\eta$ is the connection 1-form of the circle
bundle and $g_Y$ is any orbifold metric, we can form the metric
$g_X=\eta \ot \eta + \pi^*(g_Y)$ on $X$. After perturbing the
Seiberg-Witten equations on $Y$ by a closed orbifold 2-form
$\delta$ and on $X$ by its self-dual pullback $\pi^*(\delta)^+ =
\frac12(1+\star)\pi^*(\delta)$, there is a moduli space of
irreducible solutions to the Seiberg-Witten equations
$\mathcal{M}^*(X,g_X,\pi^*(\delta)^+)$ associated with $X$ and
$\mathcal{M}^*(Y,g_Y,\delta)$ associated with $Y$ (see subsections
\ref{sec:sw_eq_3_orb} and \ref{sec:sw_eq_4_man} for
definitions).  Let $\mathcal{N}^*(X,g_X,\pi^*(\delta)^+)$ be the
subcomponent of $\mathcal{M}^*(X,g_X,\pi^*(\delta)^+)$ which are
the \Spinc\ structures that are pulled back from $Y$. Theorem A
tells us that these are the only \Spinc\ structures that are
useful to study.  We can now state the main theorem of this paper:

\begin{MainThm} The pullback map $\pi^*$ induces a homeomorphism
\[ \pi^*:\mathcal{M}^*(Y,g_Y,\delta) \ra \mathcal{N}^*(X,g_X,\pi^*(\delta)^+).\]
Furthermore, if either of the two moduli spaces is a smooth
manifold, then both of them are smooth, and $\pi^*$ is a
diffeomorphism.
\end{MainThm}

The approach to the proof of Theorem B was inspired by similar
work done in \cite{sw:sw_inv_seifert_space}.

As in the free case, a manifold with a fixed point free
$S^1$-action can still be considered as a unit circle bundle, but
now it is a unit circle bundle of an orbifold line bundle over a
$3$-orbifold.  In this setup, $H^2(Y;\BZ)$ is replaced by a group
called $\Pic^t(Y)$ which records local data around the singular
set (see subsection \ref{sec:orb_line_bdle}).  Our main results
express the Seiberg-Witten invariants of $X$ in terms of the
Seiberg-Witten invariants of the orbifold $Y$ and the orbifold
Euler class $\chi$:

\begin{MainThm}
Let $X$ be a closed smooth $4$-manifold with $b_+{>}1$ and a fixed
point free circle action.  Let $Y^3$ be the orbifold quotient
space and suppose that $\chi \in \Pic^t(Y)$ is the orbifold Euler
class of the circle action. If $\xi$ is a \Spinc\ structure over
$X$ with $SW^4_X(\xi)\not=0$, then $\xi=\pi^*(\xi_0)$ for some
\Spinc\ structure on $Y$ and
\begin{eqnarray*}
SW^4_X(\xi) = \sum_{\xi' \equiv \xi_0 \mod \chi} SW^3_Y(\xi'),
\end{eqnarray*}
where $\xi'-\xi_0$ is a well-defined element of $\Pic^t(Y)$. When
$b_+{=}1$, the  formula holds for all \Spinc\ structures pulled
back from $Y$.
\end{MainThm}

This results produces two immediate corollaries.   One is a
corresponding formula for manifolds with free circle actions.
This corollary is useful for calculating examples. The second
corollary is a proof of the well known fact that the
Seiberg-Witten invariants of $Y^3\times S^1$ are the same as the
Seiberg-Witten invariants of $Y^3$.

Theorem C, together with the conjectured formula when $X$
contains fixed points (see section \ref{sec:final_remarks}), would
completely calculate the Seiberg-Witten invariants for all
$b_+>1$ 4-manifolds with circle actions. These calculations are
useful beyond just distinguishing manifolds. When Seiberg-Witten
invariants are combined with C. Taubes's results on symplectic
manifolds (c.f. \cite{sw:sw_inv_and_symp_form}), the formulas
become an easy and powerful way of calculating an obstruction for
an $S^1$-manifold to admit a symplectic structure.

We use Theorem B to study the moduli spaces in the case when
$b_+(X)=1$.  Normally when $b_+(X)=1$, the Seiberg-Witten
invariant depends on the ``chamber'' of the metric used to
calculate it.  A theorem of T. J. Li and A. Liu
\cite{sw:wallcross} shows how the numerical invariant changes when
the metric moves from one chamber into another.  Under certain
conditions, their theorem says that the invariant does not change
(making it a diffeomorphism invariant again). We show how to
construct an infinite number of $b_+=1$ manifolds with this
property. Theorem B provides a way to see explicitly why the
invariants do not change when a chamber wall is crossed.

\bigskip

\noindent{\em Acknowledgements:}\ I am deeply grateful to Ronald
Fintushel for introducing me to this field and for his constant
encouragement.  I would also like to thank Thomas H. Parker and
Peter Ozsv\'{a}th  for their many helpful discussions.

\section{Seiberg-Witten on $3$-orbifolds}

We show that all of the usual notions of gauge theory hold for
$3$-dimensional real orbifolds. Throughout, we assume that all
orbifolds are oriented, connected, and closed unless otherwise
specified.  We start with the definition of orbifolds (c.f.
\cite{orb:gauss_bonnet}).

\subsection{Definitions}

An $n$-{\em dimensional orbifold} $Y$ is a Hausdorff  space $|Y|$
together with a system $\Xi = (\{U_i \}, \{\varphi_i\},
\{\tilde{U}_i \}, \{G_i\}, \{\tilde{\varphi}_{ij} \})$ which
satisfies
\begin{enumerate}
\item $\{U_i\}$ is locally finite.
\item $\{ U_i\}$ is closed under finite intersections.
\item For each $U_i$, there exist a finite group $G_i$ acting
smoothly and effectively on a connected open subset $\tilde{U}_i$
of $\BR^n$ and a homeomorphism $\varphi_i:\tilde{U}_i/G_i \ra U_i$.
\item If $U_i \subset U_j$, there exist a monomorphism $f_{ij}: G_i \ra G_j$
and a smooth embedding $\tilde{\varphi}_{ij}: \tilde{U}_i \ra
\tilde{U}_j$ such that for all $g\in G_i$, $x \in \tilde{U}_i$,
$\tilde{\varphi}_{ij}(g \cdot x) = f_{ij}(g)\cdot
\tilde{\varphi}_{ij}(x)$ making the following diagram commute:
\[
\xymatrix{\tilde{U}_i \ar[r]^{\tilde{\varphi}_{ij}} \ar[d]_{r_i} & \tilde{U}_j \ar[d]^{r_j} \\
\tilde{U}_i / G_i \ar[r]^{\varphi_{ij}} \ar[d]_{\varphi_i} & \tilde{U}_j /G_j \ar[d]^{\varphi_j}\\
U_i \ar[r] & U_j
}
\]
where  $\varphi_{ij}$ are induced by the monomorphisms and the
$r_i$'s are the natural projections.
\end{enumerate}

The system $\Xi$ is called an {\em atlas} and each $\varphi_i
\circ r_i:\tilde{U}_i \ra U_i$ is called a {\em local chart}.  An
orbifold $Y$ is {\em connected}  and {\em closed} if the
underlying space $|Y|$ is.  Two atlases give the same orbifold
structure if there is a common refinement.

Let $x\in |Y|$ and $\tilde{U}_x \ra U$ be a local chart
containing $x$.  The local group at $x$, denoted $G_x$, is the
isotropy group of $G$ of any point in $U$ corresponding to $x$
(well-defined up to isomorphism).  Set $\Sigma Y = \{ x\in |Y|
\;\; | \;\; G_x \not= 1 \}$.  This set is closed and nowhere
dense, and in fact it is easily shown that $\dim \Sigma Y \leq
n-2$. After removing the singular set, $Y \setminus \Sigma Y$
becomes a manifold.

All theorems henceforth will be stated and proved for
$3$-dimensional orbifolds $Y$  where $\Sigma Y$ is a  finite
disjoint set of smooth closed curves $l_1,\dots, l_n$ that are
assigned integral multiplicities $\alpha_1, \dots, \alpha_n$
given by their local isotropy group $\BZ_{\alpha_i} =
\BZ/\alpha_i\BZ$. Let $D$ be the standard complex disk and
consider a $\BZ_{\alpha_i}$ action on it by rotation.  We will
take a convenient atlas in all of the atlases which give the same
orbifold structure. Equip $|Y|$ with an atlas of coordinate charts

\begin{eqnarray*}
\phi_i: (S^1 \times D, S^1 \times 0) \ra (U_i, l_i) & i = 1, \dots, n \\
\phi_x: D^3_x \ra U_x & x \in  Y \setminus \{ l_1, \dots, l_n\}, \\
\end{eqnarray*}
where the $\phi_i$ induce homeomorphisms from $(S^1 \times
D/\BZ_{\alpha_i}, S^1 \times 0)$ to $(U_i, l_i)$, the $\phi_x$
are homeomorphisms, the $U_i$ are all pairwise disjoint, $U_x
\cap \Sigma Y = \emptyset$, and the transition functions are all
diffeomorphisms.

\begin{example}
The triple $Y=(S^3, K, n)$ where $K$ is a knot in $S^3$, $K$ is
the singular locus $\Sigma Y = K$, and the isotropy group around
$K$ is $\BZ_n$, is an example of a $3$-orbifold.
\end{example}

Define an $n$-dimensional {\em orbifold bundle} over $Y$ in the
following manner.  Set $U_x \times V^n$ over each  $U_x$ for an
$n$-dimensional vector space $V^n$.  Over $U_i$ the vector bundle
is given by the quotient $(S^1 \times D \times
V^n)/\BZ_{\alpha_i}$ where $(S^1\times D \times V^n)$ is a
$\BZ_{\alpha_i}$-equivariant vector bundle specified up to
isometry by giving a representation $\sigma_i: \BZ_{\alpha_i} \ra
GL_n(V)$. The vector bundle over $Y$ is then specified by a
1-cocycle of transition functions over the overlaps.

\subsection{Orbifold line bundles}
\label{sec:orb_line_bdle}

Under tensor product the topological isomorphism classes of
orbifold line bundles form a group $\Pic^t(Y)$ called the {\em
topological Picard group}.  We describe this group in this
subsection.

We can record the information in $\Pic^t(Y)$ by using a
generalization of equivariant cohomology. Think of  $Y$ as the
union of $Y\setminus \{l_1,\dots,l_n\}$ and $\displaystyle \amalg
\; (l_i \times D/\BZ_{\alpha_i}) $. Define $Y_V$ to be the union
of $Y\setminus \{l_1,\dots, l_n\}$ and $\displaystyle \amalg \;
\left(l_i \times (D\times_{\BZ_{\alpha_i}}
E\BZ_{\alpha_i})\right)$ glued using sections of $l_i \times
(D\setminus \{0\} \times_{\BZ_{\alpha_i}} E\BZ_{\alpha_i}) \ra
U_i \setminus l_i$.  These sections are unique up to homotopy
because the fibers of the bundle are contractible.

The following theorem is contained in
\cite{circact:sf_hom_3_sphere}.

\begin{thm} The following groups are isomorphic:
\begin{enumerate}
\item $H^1(Y_V;\BZ) \cong H^1(|Y|;\BZ)$,
\item $H^2(Y_V;\BZ) \cong \Pic^t(Y)$.
\end{enumerate}
\label{thm:isom_cohom}
\end{thm}

\begin{rem}
In the literature, the group $H^*_V(Y) := H^*(Y_V)$ is often
called the {\em V-cohomology ring of $Y$}.
\end{rem}

Here is another way to describe $\Pic^t(Y)$.  For a fixed $i$,
define an orbifold line bundle over $Y$ to be a trivial line
bundle $A = (Y\setminus l_i) \times \BC$ and over  $U_i$ it is
given by $B=(S^1 \times D \times \BC_{\xi})/\BZ_{\alpha_i}$ where
$a \in \BZ_{\alpha_i}$ acts using the standard representation
\[ a\cdot (\gamma, w,z) \mapsto (\gamma, \euler^{\frac{2\pi \I
a}{\alpha_i}}w,\; \euler^{\frac{2\pi \I a}{\alpha_i}} z).
\]
The bundle is glued together using a transition function
$\varphi_{BA}(\gamma,w) = w$ on the overlap $S^1\times
(D\setminus \{0\})$.  For each $l_i$, create such a line bundle
called $E_i$.

Let $L$ be an orbifold line bundle over $Y$.  There is a
collection of integers $\beta_1,\ldots\beta_n$  satisfying $$0
\leq \beta_i < \alpha_i$$ such that the bundle $L\ot
E_1^{-\beta_1} \ot \cdots \ot E_n^{-\beta_n}$ is a trivial
orbifold line bundle over each neighborhood of the $l_i$'s. By
forgetting the orbifold structure, it can be naturally identified
with a smooth line bundle (denoted by $|L|$) over the smooth
manifold $|Y|$.

\begin{thm}
The isomorphism classes of orbifold line bundles on $Y$ with
specified isotropy  representations $\xi_1^{\beta_1}, \ldots,
\xi_n^{\beta_n}$ along $l_1,\ldots, l_n$ respectively are in
bijective correspondence with $\chi \in H^2(|Y|;\BZ)$.
\label{thm:orb_iso_class}
\end{thm}

The proof below generalizes~\cite{circact:Class_of_circ_4man} to
the case of an arbitrary orbifold line bundle.

\begin{proof}
Given $L \in \Pic^t(Y)$, we construct $L\ot E_1^{-\beta_1}
\ot\cdots \ot E_n^{-\beta_n}$ and its desingularization $|L|$
explicitly. Let $\pi:X \ra Y$ be the unit circle bundle of $L$.
Set $Q = \cup \, U_i$ in $Y$ and $P= \pi^{-1}(Q)$ with
$P_i=\pi^{-1}(U_i)$. Then $X' = X\setminus P$ is a principal
$S^1$-bundle over $Y' = Y\setminus Q$.

In general, the unit circle bundle $X$ is an orbifold rather than
a manifold.  $P_i$ is a quotient of $ l_i \times D^2\times S^1$
by the action of $\BZ_\alpha$ defined by $\xi:(\gamma, w, t)
\mapsto (\gamma, \xi w, \xi^\beta t )$. It follows that the
isotropy group of a point in the quotient of $l_i \times
\{0\}\times S^1$ is $\{\xi\in \BZ_\alpha:\xi^\beta = 1\}$ for all
points $p\in l$. When the isotropy group is trivial
($gcd(\alpha,\beta)=1$) the quotient is smooth. In the case that
$\beta \equiv 0 \mbox{ mod}_\alpha$, $L$ is a usual line bundle
around that loop, but the $4$-manifold still has a nontrivial
orbifold structure.  Set  $d= gcd(\alpha_i,\beta_i)$.

Let $m_i = \partial (\{0\} \times D \times \{1\})$  be the
meridian loop of $l_i$ before the quotient is taken. Denote the
class it represents in the quotient $P_i$ by $\tilde{m}_i$.  The
homeomorphism $\varphi_i:\partial P_i \ra \partial X'$ determines
a section $s:\partial Y' \ra \partial X_0$ which is specified up
to homology by the relation:
$$\varphi_*[\tilde{m}_i] = \left(\frac{\alpha_i}{d}\right) s_*[m_i'] + \left(\frac{\beta_i}{d}\right) [f'],$$
where $m_i'$ is the meridian of $l_i$ in $Y'$,  $f'$ is a fiber
of $\partial X'$, and $0\leq \beta_i \leq \alpha_i$.    The local
invariants $(\alpha_i,\beta_i)$ specify $P_i$ up to
orientation-preserving equivariant homeomorphism.

The bundle $X'$ can be extended to the unit circle bundle of $|L|$
by equivariantly attaching $l_i\times D \times S^1$ with a bundle
isomorphism $\phi_i$. Bundle isomorphisms covering the identity
are classified up to vertical equivariant isotopy by homotopy
classes of maps in $[\partial(S^1 \times D),S^1]=\BZ\oplus\BZ$.
However we can change $\phi_i$ by a bundle automorphism
classified by $[S^1\times D,S^1]=H^1(S^1\times D;\BZ)=\BZ$; these
maps  change $(\phi_i)_*([l_i])$ by a multiple of the fiber.
Therefore the resulting bundle $X'\cup_{\phi_i} (l_i\times
D\times S^1)$ can be completely specified by the map
$$(\phi_i)_*[m_i]=s_*[m']+r[f']$$ for some $r\in \BZ$.  Thus we
determine the principal $S^1$ bundle of $|L|$ by specifying that
$r=0$.

In summary:

\begin{enumerate}
\item The unit circle bundle of $L$ is obtained by gluing the
quotient $P$ using maps
\[(\varphi_i)_*[\tilde{m}_i] = \left(\frac{\alpha_i}{d}\right)
s_*[m_i'] + \left(\frac{\beta_i}{d}\right) [f'].\]  Note that this
bundle depends only on the section $s_*[m_i']$ as well.


\item The unit circle bundle of $L\ot E_1^{-\beta_1}\ot \cdots \ot
E_n^{-\beta_n}$ is obtained by gluing in the quotient $\amalg_i \;
l_i\times (D/\BZ_{\alpha_i}) \times S^1$ into $X'$ using  maps
$(\phi_i)_*[\tilde{m}_i] = s_*[m_i'].$

\item The unit circle bundle of the desingularization $|L|$ is
obtained by gluing in $\amalg_i \; l_i\times D\times S^1$ using
maps $(\phi_i)_*[m_i] = s_*[m_i'].$
\end{enumerate}

Next we show that two orbifold line bundles $L_1$ and $L_2$ with
the same isotropy representations and equivalent
desingularizations $|L_1|= |L_2|$ are equivalent as orbifold line
bundles.

Construct two principal $S^1$-bundles $X_1$ and $X_2$ from $X'$
to form unit circle bundles $|L_1|$ and $|L_2|$. The construction
depends on choices of the class $\sum_{i=1}^n (s_j)_*[m_i'] \in
H_1(\partial X';\BZ)$ coming from sections $s_j:\partial Y' \ra
\partial X_j'$ for $j=1,2$.  Let $\theta_j \in H^2(Y',\partial
Y')$ be the obstruction to extending these sections over $X_j'$.
Let $\tau \in H^1(\partial Y')$ be the primary difference of
$s_1$ and $s_2$.  A diagram chase
\[\xymatrix{
H^1(\amalg_i l_i\times D;\BZ) \ar[r]^{\delta_1} \ar[d]_{i^*} &
H^2(Y,\amalg_i l_i\times D;\BZ) \ar[r]^{\hspace{.8cm} j_1^*}
\ar[d]_{\lambda}^{\cong} &
H^2(Y;\BZ) \\
H^1(\partial Y';\BZ) \ar[r]^{\delta_2} \ar[d]_{\delta_3} &
H^2(Y', \partial Y';\BZ) \\
H^2(\amalg_i l_i \times D, \partial;\BZ)}\] shows that
$j_1^*\lambda^{-1}\delta_2(\tau) =
j_1^*\lambda^{-1}(\theta_1-\theta_2)= c_1(|L_1|) - c_1(|L_2|)=0$.
Thus there is an element $\tau' \in H^1(\amalg_i l_i\times D)$
such that $\delta_1\tau' = \lambda^{-1}\delta_2\tau$, and
$\delta_3(\tau - i^*\tau')=0$. Therefore $\tau \in
i^*(H^1(\amalg_i l_i\times D;\BZ))$ implying that $(s_1)_*[m_i']$
is homotopic to $(s_2)_*[m_i']$ through a homotopy in $l_i \times
D$.  Since the construction of the unit circle bundle of the
orbifold line bundle in (1) depended only on these sections,
$L_1$ and $L_2$ are equivalent.
\end{proof}

The above theorem means that a given orbifold line bundle $L$
over $Y$ is specified by the data
$$ (c_1(|L|), \beta_1, \ldots, \beta_n)$$ called the {\em Seifert
invariant of $L$ over $Y$}.  (This data, of course, is not
unique).

\subsection{Spin$^c$ Structures on 3-orbifolds}

The \Spinc\ structures on a 3-orbifold $Y$ are defined by a pair
$\xi = (W, \rho)$ consisting of a rank 2 complex orbifold bundle
$W$ with a hermitian metric (the spinor bundle) and an action
$\rho$ of orbifold 1-forms on spinors,
\[\rho: T^*Y \ra \mbox{End}(W),\]
which satisfies the property that, if $e^1, e^2, e^3$ are an
orthonormal coframe at a point in $Y$, then the endomorphisms
$\rho(e^i)$ are skew-adjoint and satisfy the Clifford relations
\[\rho(e^i)\rho(e^j) + \rho(e^j)\rho(e^i) = -2\delta_{ij}.\]
We also require that the volume form $e^1 \wedge e^2 \wedge e^3$ acts by
\[\rho(e^1\wedge e^2 \wedge e^3) = -\Id_{W}.\]
We will write $c_1(\xi)$ for the first Chern class of $\det W$.

\begin{thm}
The tangent bundle $T^*Y$ of an orbifold always lifts to an
orbifold $\Spinc(3)$-bundle. \label{thm:lift_to_spinc}
\end{thm}

\begin{proof}
If we can split $TY$ into a 1-dimensional real line bundle and a
complex orbifold line bundle $L$, then $w_2(TY) = w_2(\BR \oplus
L) = w_2(L)$ is just the mod 2 reduction of $c_1(L)$ for some
orbifold line bundle $L$.  Hence $TY$ lifts.

We need to find a nowhere zero section of $TY$.  Note that each
$l_i\times D/\BZ_{\alpha_i}$ comes with an
$\BZ_{\alpha_i}$-invariant oriented nonzero vector field that is
tangent to $l_i$ at each point in $D/\BZ$. This vector field
induces a nonzero section $s:\partial Y' \ra TY|_\partial$.
Remove an extra $S^1\times D$ from the interior of $Y'$ and put
 a similar nonzero section on the boundary. The obstruction to
extending the section into the interior of
$$Y'' = Y\setminus \left((S^1\times D) \;\; \cup \;\; \amalg_i \;
l_i\times D/\BZ_{\alpha_i}\right)$$ is an element of
$H^3(Y'',\partial Y'';\pi_2(S^2)) = \BZ.$ Using the homology
relation
$$[\partial (S^1\times D)] = - \sum_i [\partial(l_i\times D/\BZ_{\alpha_i})],$$
the obstruction can be removed by changing the framing on the
boundary of $S^1\times D$.  Thus $TY$ admits a nowhere zero
vector field.
\end{proof}

\begin{rem} In \cite{orb:gauss_bonnet},
I. Satake treated the V-Euler class as the index of a unit vector
field on $TY$ with singularities and showed that $\chi_V(Y)=0$ for
odd dimensional orbifolds.  Thus it is not surprising that nonzero
vector fields exists on 3-orbifolds.
\end{rem}

\begin{thm} The set of \Spinc\ structures lifting the frame bundle
of a $3$-orbifold $Y$ is a principal homogeneous space over
$\Pic^t(Y)$: The difference of two \Spinc\ structures
$\xi_1,\xi_2$ is an orbifold line bundle.
\end{thm}

\begin{proof}
Let $\xi_1$ and $\xi_2$ be two \Spinc\ structures which are lifts
of the frame bundle.  Away from the $l_i$'s, the difference of
two \Spinc\ structures is a complex line bundle as in the smooth
case.   Because
$$\left(c_1(\xi_1)-c_1(\xi_2)\right)[\partial(l_i \times D/\BZ_{\alpha_i})] = 0$$ for
all $l_i$, we can extend the complex line bundle over the
desingularization $|Y|$ using techniques in Theorem
\ref{thm:orb_iso_class}. Thus we can investigate locally to show
that any two lifts of isotropy representations into $\Spinc(3)$
differ by a representation into $S^1$.  Note that this is not
immediately obvious because there are many different
representations of $\BZ_\alpha$ into $\Spinc(3)=U(2)$.

Let $\Theta_i$ be the unit vector field on $l_i\times
D/\BZ_{\alpha}$ which is tangent to the circle $l_i$ at each
point. We use the fact that $\rho:\BZ_{\alpha} \ra SO(3)$ is a
rotation which leaves the nonzero vector field $\Theta_i$
invariant. Identify $SU(2)$ with the unit quaternions. The map
$Ad:SU(2) \ra SO(3)$ is given by
$$g \mapsto gh\bar{g}$$ for all $h \in \Imaginary \BH$ and is the double
cover of $SO(3)$.  Thus $SO(3)$ can be thought of as the unit
quaternions modulo the equivalence $h\sim -h$. Without loss of
generality, we may assume that the invariant vector field
$\Theta_i$ is generated by ${\bf i} \in \Imaginary \BH$ at each
point in $l_i$.

It is easy to see that elements of $SO(3)$ which rotate the
second two components while leaving $\I$ invariant are of the
form $e^{\I\theta} \in \BH$.  Hence $\rho(1)=\lambda^{\tau}$ where
$\lambda$ is a $2\alpha$-root of unity in $\BC$ and $0\leq \tau <
\alpha$.

The \Spinc\ representation $\sigma:\Spinc \ra \End(\BH)$ is given
by $\sigma(g,e^{i\theta})h = ghe^{i\theta}$ for all $h\in \BH$.
Here we have used the fact $\Spinc(3) = SU(2) \times S^1 \, /\,
(-1,-1)$. Using this identification, \Spinc\ projects to $SO(3)$
by the adjoint map as well,
$$(g,e^{i\theta}) \mapsto gh\bar{g}$$ for all $h \in \Imaginary \BH$.
Thus the representation $\rho$ lifts to $\tilde{\rho}$:
\[\xymatrix{
& S^1 \ar[d] \\
& \Spinc\ \ar[d]\\
\BZ_{\alpha} \ar[r]_{\rho} \ar[ru]_{\tilde{\rho}}
\ar[ruu]^{\hat{\rho}} & SO(3)}\] given by $\tilde{\rho}(1) =
(\lambda^{\tau}, \hat{\rho}(1))$ \ (or equivalently
$(-\lambda^{\tau},-\hat{\rho}(1))$) for some representation
$\hat{\rho}:\BZ_{\alpha} \ra S^1$.  The representation
$\hat{\rho}$ is given by $\hat{\rho}(1)=\lambda^\kappa$ for some
$0\leq \kappa <\alpha$. Hence the difference of two  \Spinc\
structures $\xi_1- \xi_2$ locally is a representation
$(\hat{\rho}_1 - \hat{\rho}_2):\BZ_{\alpha_i} \ra S^1$.

Globally, $\xi_1-\xi_2$ differs by a complex line bundle over
$|Y|$ and local isotropy representations into $S^1$, i.e., an
element in $\Pic^t(Y)$ as described in Theorem
\ref{thm:orb_iso_class}.
\end{proof}

\subsection{Seiberg-Witten Equations on 3-orbifolds}
\label{sec:sw_eq_3_orb}

Fix an orbifold $SO(3)$-connection $\Snabla$ on the cotangent
bundle $T^*Y$ and a \Spinc\ structure $\xi =(W,\rho)$.

\begin{df}
A Hermitian connection $\Spnabla$ on $W$ is called {\em spinorial}
with respect to $\Snabla$  if it is compatible with Clifford
multiplication, i.e.,
\begin{eqnarray}
\Spnabla(\rho(v)\psi) = \rho(\Snabla v)\psi + \rho(v)(\Spnabla
\psi). \label{eq:spinorial}
\end{eqnarray}
The set of all spinorial connections will be denoted by $\mathcal{A}(W)$.
\end{df}

Given a trivialization for $W$, the connection matrix of any
$\Snabla$-spinorial connection $\Spnabla$ can be written with
respect to this trivialization as
\[\label{eq:conn_matrix} \frac14 \sum \omega^i_j \ot \rho(e^i\wedge e^j) +  b\Id_W\]
where $\omega^i_j$ are the connection matrices for $\Snabla$, and
$b\in \Omega^1(Y,\I\BR)$ is an orbifold 1-form.   We will often
think of spinorial connections on a \Spinc\ structure as $U(1)$
connections on det $W$ coupled with the Levi-Civita connection
$\nabla_Y$ on $T^*Y$.  A spinorial connection $\Spnabla$ defines a
Dirac operator $D_A: \Gamma(Y,W) \ra \Gamma(Y,W)$ on the space of
orbifold sections of $W$ which is self-adjoint. The perturbed
Seiberg-Witten equations are the following pair of equations for
$(A,\Psi)$ where $A$ is a $U(1)$ orbifold connection on det $W$
and $\Psi$ is an orbifold section of $W$:

\begin{eqnarray}
\label{eq:sw_3dim}
\begin{array}{r@{=}l}F_A + \delta  -\star \tau(\Psi) \ \ & \ \ 0\\
D_A(\Psi) \ & \ 0.
\end{array}
\end{eqnarray}
Here $\tau:\Gamma(Y,W) \ra \Omega^1(Y;\I\BR)$ is the adjoint to
Clifford multiplication, defined by
\[ \langle \rho(b) \Psi,\Psi\rangle_W = 2\langle b,\tau(\Psi) \rangle_{\Lambda^1},\]
for all orbifold $1$-forms $b\in \Omega^2(Y;\I\BR)$ and all $\Psi
\in \Gamma(Y,W)$. The $\delta\in \Omega^2(Y;\I\BR)$ is a closed
orbifold 2-form used to perturb the equations.

For a fixed metric $g_Y$ and perturbation term $\delta$, the
moduli space $\mathcal{M}(Y,\xi, g_Y,\delta)$ is the space of
solutions to (\ref{eq:sw_3dim}) modulo the action of the gauge
group $\mathcal{G}=Map(Y,S^1)$.    Let
$\mathcal{M}^*(Y,\xi,g_Y,\delta)$ denote the set of irreducible
solutions (i.e., where $\Psi \not \equiv 0$). For a generic
perturbation, the moduli space is a compact, smooth manifold
containing no reducible solutions.  In that case, the fundamental
class $[\mathcal{M}(Y,\xi, g_Y, \delta)]$ is essentially the
Seiberg-Witten invariant. Evaluating it against some universal
classes defines a map $SW^3(\xi) \in \BZ$ which is independent of
the Riemannian metric and perturbation when $b_1(Y)>1$ (c.f.,
~\cite{sw:sw_eq_and_fourmanifolds}). Denote the union over all
distinct \Spinc\ structures by $\mathcal{M}(Y,g_Y,\delta)$.

\section{$4$-Manifolds with fixed point free circle actions}

In this section we study manifolds with fixed point free circle
actions.  We describe the cohomology of these manifolds and show
that under some circumstances, line bundles with connections can
be pushed forward to orbifold line bundles with connection on the
quotient.  Finally, we describe how to pullback \Spinc\ structures
from \Spinc\ structures on the quotient.

\subsection{Homology}
A $4$-manifold with fixed point free $S^1$-action can be viewed
as the boundary of a disk bundle or the unit circle bundle of an
orbifold line bundle $L$ over a $3$-orbifold $Y$.   Henceforth, we
will assume that $X$ is a unit circle orbifold line bundle $L$
over $Y$ where each local invariant $\beta_i$ is relatively prime
to $\alpha_i$. Denote $\pi:X \ra Y$ for the projection map.

When $X$ is smooth, then $X_V \ra Y_V$ is an honest $S^1$-bundle
and we have the Gysin sequence:
\[
\xymatrix{
0 \ar[r] & H^1_V(Y) \ar[r] \ar@{=}[d] & H^1_V(X) \ar[r] \ar@{=}[d] & H^0_V(Y) \ar[r] \ar@{=}[d] & H^2_V(Y) \ar@{=}[d] \\
0 \ar[r] & H^1(|Y|) \ar[r] & H^1(X) \ar[r] & \BZ \ar[r] & \Pic^t(Y)\\
& & & 1 \ar[r] & [L] \\
\ar[r] & H^2_V(X) \ar[r] \ar@{=}[d] & H^1_V(Y) \ar[r] \ar@{=}[d] & H^3_V(Y) \ar@{=}[d] \\
\ar[r] & H^2(X) \ar[r] & H^1(|Y|) \ar[r] & H^3_V(Y) }
\]

\begin{thm}\label{thm:X_cohomology}
If $X$ is a $4$-manifold with a fixed-point free circle action
over $Y$ given by the sphere bundle of a line bundle $L$ over
$Y$, then
\begin{eqnarray*}
H^1(X,\BZ) & \cong & \left\{ \begin{array}{rl} H^1(|Y|;\BZ), & [L] \mbox{ is not torsion }\\ H^1(|Y|,\BZ) \oplus \BZ, & [L] 
\mbox{ is torsion}\end{array}
\right. \\
H^2(X;\BZ) &\cong  & (\Pic^t(Y)/<[L]>) \oplus \ker (\cdot \cup
[L]):H^1(|Y|;\BZ)\ra H^3_V(Y;\BZ).
\end{eqnarray*}
In particular, since the kernel of $(\cdot \cup [L]):H^1\ra H^3$
is torsion free, all torsion classes must come from pullbacks in
$\pi^*(\Pic^t(Y))$.
\end{thm}

When $[L]$ is not torsion, the rank of $\Pic^t(Y)/<[L]>$ and
$\ker (\cdot \cup [L])$ are both equal to $b_1(|Y|)-1$.  A basis
for the former space can be represented by the Poincar\'{e} duals
of tori of the form $\pi^{-1}(loop)$ for smooth loops in
$Y\setminus \Sigma Y$.  A basis for the later space can be
represented by surfaces in $X$ which, after integrating over the
fiber, are the Poincar\'{e} duals of surfaces in $|Y|$.  The
simple intersection relationship between loops and surfaces in
$|Y|$ implies that the intersection form $Q_X$ should be simple
as well.

In fact, since the signature is zero (c.f.
\cite{circact:circ_act_on_4man_II}), the classification of
intersection forms says that $Q_X$ is equivalent to the direct sum
of matrices of the form (where $d$ an integer)
$$\displaystyle \left(\begin{array}{cc} 0 & 1\\ 1& d\end{array}\right)$$
with respect to a basis $\{A,B\}$ where $A\in \pi^*(\Pic^t(Y))$ is
a class pulled back from $Y$.  Pulled back classes always have
square zero by the naturality of the cup product and the fact
that the product of 2-forms on $Y$ is always zero.

\subsection{Line bundles over $X$}
\label{sec:linebundlesoverX}

Orbifold line bundles $E$ over $Y$ pullback to usual line bundles
$\pi^*(E)$ over $X$.  Except for the case $X = |Y| \times S^1$,
this is a many to one correspondence.  Nonetheless, it can be made
faithful in the following way.  Given a line bundle $E$ with
connection $A$ over $X$ with the following two properties:
\begin{enumerate}
\item The curvature two form of $A$ pulls up from $Y$, i.e.,
$$\iota_{\T} F_A = 0,$$ where
$\T$ is the everywhere non-zero vector field generated by the
circle action on $X$.
\item There exists a point $x \in Y \setminus \Sigma Y$ such
that holonomy of $A$ around $\pi^{-1}(x)$ is trivial.
\end{enumerate}
Then $(E,A)$ can be pushed forward to an orbifold line bundle
with connection on $Y$ (up to gauge equivalence).    If one such
point $x \in Y\setminus \Sigma Y$ satisfies the second condition,
then all points outside the critical set do.  Such connections
are said to have {\em trivial fiberwise holonomy}.

We state Proposition 5.1.3 from~\cite{sw:sw_inv_seifert_space}.

\begin{prop} There is a natural one-to-one correspondence between
orbifold line bundles with connection over $Y$ and usual line
bundles with connection over $X$, whose curvature forms pull up
from $Y$ and whose fiberwise holonomy is trivial.  Furthermore,
this correspondence induces an identification between orbifold
sections of the orbifold bundle over $Y$ with fiberwise constant
sections of its pullback over $X$. \label{prop:1to1}
\end{prop}

Pull back connections $\pi^*A$ are  characterized by
$\Spnabla^{\pi^*A}_{\T}\Psi = 0$ for all pulled back sections
$\Psi$.

\subsection{Seiberg-Witten Equations of Smooth 4-manifolds}
\label{sec:sw_eq_4_man}

A \Spinc\ structure $\xi=(W,\sigma)$ on an oriented $4$-manifold
$X$ is a hermitian vector bundle $W$ of rank $4$, together with a
Clifford multiplication $\sigma:T^*X\ra End(W)$.  The bundle $W$
decomposes into two bundles of rank $2$, $W^+ \oplus W^-$, with
$\det W^+ = \det W^-$. The bundle $W^-$ is the subspace
annihilated by the action of self-dual $2$-forms. We set
$c_1(\xi)$ to be the first Chern class of $\det W^+$.

There is a natural way to pullback a $\Spinc$ structure from $Y$
to $X$. Let $\I\eta$ denote the connection 1-form of the circle
bundle $\pi:X \ra Y$, and let $g_Y$ be a metric on $Y$, then endow
$X$ with the metric $g_X = \eta \ot \eta + \pi^*(g_Y)$.  Using
this metric, there is an orthogonal splitting
\[ T^*X \cong \BR\eta \oplus \pi^*(T^*Y).\]

If $\xi = (W,\rho)$ is a \Spinc\ structure over $Y$, define the pullback
of $\xi$ to be $\pi^*(\xi) = (\pi^*(W) \oplus \pi^*(W), \sigma)$ where the action
\[ \sigma:T^*X \ra \mbox{End}(\pi^*(W)\oplus \pi^*(W))\]
is given by
\[\sigma(b\eta + \pi^*(a)) = \left( \begin{array}{cc} 0 & \pi^*(\rho(a)) +
b\Id_{\pi^*(W)} \\ \pi^*(\rho(a)) - b\Id_{\pi^*(W)} & 0 \end{array} \right) .\]
This defines a \Spinc\ structure on $X$.

Choosing a \Spinc\ structure $\xi_0 = (W_0,\rho)$ on $Y$ gives
rise to a one-to-one correspondence between Hermitian orbifold
line bundles and \Spinc\ structures on $Y$ via $E \mapsto W_0 \ot
E$. Likewise, the pullback \Spinc\ structure $\xi=\pi^*(\xi_0)$
induces a  one-to-one correspondence between Hermitian line
bundles and \Spinc\ structures on $X$.

\begin{rem}
\label{rem:push_spinc} In this way we can think of a \Spinc\
structure with respect to $\xi_0$ or $\xi$ as a choice of line
bundle on $Y$ or $X$ respectively.  This allows us to push-forward
a \Spinc\ structure with a trivial fiberwise connection on $\det
W^+$ from $X$ to $Y$ via Proposition~\ref{prop:1to1}.
\end{rem}

There is a natural connection on $X$ which is compatible with the
reduction $T^*X = \BR \eta \oplus \pi^*(T^*Y$).  Let $\nabla^Y$
denote the Levi-Civita connection on $Y$ and set $\Snabla = d
\oplus \pi^*(\nabla^Y)$.  This is a compatible connection which
satisfies
\begin{eqnarray}
\Snabla\eta = 0, \ \ \  \mbox{ and } \ \ \ \Snabla(\pi^*(\beta))=
\pi^*(\nabla^Y\beta). \label{eq:_snabladef}
\end{eqnarray}

It is more convenient  to use this reducible $SO(4)$-connection
instead of the Levi-Civita connection. By coupling it with a
$U(1)$-connection $A$ on $\det W^+$ we can define a spinorial
connection on $W^+$.  Define a Dirac operator
$\dirac^+_A:\Gamma_X(W^+)\ra \Gamma_X(W^-)$ from the space of
smooth sections of $W^+$ to $W^-$.  The $4$-dimensional perturbed
Seiberg-Witten equations for a section $\Psi \in \Gamma_X(W^+)$
and a $U(1)$-connection $A$ on $\det W^+$ are:
\begin{eqnarray}
\label{eq:sw_4dim}
\begin{array}{r@{=}l}
F^+_A + \delta  - q(\Psi) \ \ & \ \ 0,\\
\dirac^+_A(\Psi) \ \ & \ \ 0.
\end{array}
\end{eqnarray}
Here $F^+_A$ is the projection of the curvature onto the self-dual
two forms, $\delta \in \Omega^+(X;\I\BR)$ is self-dual $2$-form
used to perturb the equations, and $q:\Gamma_X(W^+) \ra
\Omega^+(X, \I\BR)$ defined by $q(\Psi) = \Psi \ot \Psi^* -
\frac12 |\Psi|^2$ is the adjoint of Clifford multiplication by
self-dual $2$-forms,i.e,
\begin{eqnarray}
\label{eq:adjointprop}
 \langle \sigma(\beta) \Psi, \Psi\rangle_{W^+} & = & 4 \langle \beta, q(\Psi)\rangle_{\I\Lambda^+}
\end{eqnarray}
for all self-dual $2$-forms $\beta \in\Omega^+(X;\I\BR)$ and all
sections $\Psi$.

Similar to the $3$-dimensional case, the moduli space
$\mathcal{M}(X,\xi,g_X,\delta)$ is the space of solutions
$(A,\Psi)$ modulo the action of the gauge group.  We are using a
reducible connection $\Snabla$ instead of the Levi-Civita
connection on $T^*X$, but this alternative compatible connection
is an allowable perturbation of the usual Seiberg-Witten
equations and can be used to calculate the Seiberg-Witten
invariants (see section 4 of \cite{sw:symp_thom_conj}). Under
suitable generic conditions the moduli space is a compact,
oriented, smooth manifold of dimension
\begin{eqnarray}
\label{eq:dim_count} d(\xi)= \frac14\left( c_1(\xi)^2 - 2 \chi(X)
-3\sigma(X)\right)
\end{eqnarray}
which is independent of metric and perturbation when $b_+(X)>1$.

The Seiberg-Witten invariant $SW_X(\xi)$ is a suitable count of
solutions. Fix a base point in $M$ and let $\mathcal{G}^0 \subset
Map(X,S^1)$ denote the group of maps which map that point to 1.
The base moduli space, denoted by $\mathcal{M}^0$, is the
quotient of the space of solutions by $\mathcal{G}^0$. When the
moduli space $\mathcal{M}(X,\xi,g_X,\delta)$ is smooth,
$\mathcal{M}^0$ is a principle $S^1$-bundle over
$\mathcal{M}(X,\xi,g_X,\delta)$. For a given \Spinc\ structure
$\xi$, the 4-dimensional Seiberg-Witten invariant $SW_X(\xi)$ is
defined to be $0$ when $d(\xi)<0$, the sum of signed points when
$d(\xi)=0$, or if $d(\xi)>0$, it is the pairing of the
fundamental class of $\mathcal{M}(X,\xi,g_X,\delta)$ with the
maximal cup product of the Euler class of the $S^1$-bundle
$\mathcal{M}^0$.

The dimension formula (\ref{eq:dim_count}) simplifies when the
manifold has a fixed point free circle action.  Because $X$ has a
nonzero vector field $\T$, the Euler class is zero. As mentioned
previously, the signature of $X$ is also zero.

\begin{prop}\label{ref:dim_=c_1^2}
Suppose that $X$ is a 4-manifold with a fixed point free circle
action.  The expected dimension of the moduli space for a \Spinc\
structure $\xi$ is
$$d(\xi)=\frac14 c_1(\xi)^2.$$
\end{prop}

\section{Spin$^c$ structures and SW solutions}

We continue to work  with an orbifold circle bundle $\pi:X\ra Y$
with an $S^1$-invariant metric $g_X = \eta^2 +\pi^*(g_Y)$. The
perturbation $\delta \in \Omega^2(Y,\I \BR)$ is a closed orbifold
2-form used to perturb the $3$-dimensional equations which is then
pulled back and projected on to the self-dual 2-forms of $X$ to
perturb the 4-dimensional equations.

\subsection{Restrictions on Spin$^c$ structures}
\label{sec:sw_pullbacks}

First, we make some basic observations.  If $SW_X(\xi)\not=0$ for
some \Spinc\ structure $\xi$, then the expected dimension of the
moduli space  is nonnegative, implying
\begin{eqnarray}\label{eq:c_1^2geq0}
c_1(\xi)^2\geq 0. \end{eqnarray} If $b_+(X)=1$, then the metric
$g_X$ induces a splitting $\mathcal{H}^2(X;\BR) =
\mathcal{H}^+\oplus \mathcal{H}^-$ where $\mathcal{H}^+$ is one
dimensional.  Let $c_1(\xi)^+$ be the $L^2$ projection onto the
self-dual subspace $\mathcal{H}^+$.  When $c_1(\xi)^+$ is
nonzero, it provides an orientation for $\mathcal{H}^+$.  In this
situation the Seiberg-Witten invariant depends on the chamber of
$2\pi c_1(\xi)+\pi^*(\I\delta)$.  We will say that $\alpha\in
H^2(X;\BR)$ lies in the positive chamber if $\alpha^+\cdot
c_1(\xi)^+>0$.  Denote the Seiberg-Witten invariant calculated
for $\alpha = (2\pi c_1(\xi)+\pi^*(\I\delta))$ in this chamber by
$SW_X^+(\xi)$ and denote the invariant of the other chamber by
$SW_X^-(\xi)$.

When $c_1(\xi)^+=0$ there no distinguished chamber. However, if
$SW_X(\xi)\not=0$ in either chamber, the dimension of the moduli
is nonzero and
$$0\leq c_1(\xi)^2 = (c_1(\xi)^-)^2 \leq 0.$$
Since the intersection form on $\mathcal{H}^-$ is definite,
$c_1(\xi)$ is a torsion class and pulled back from $Y$ by Theorem
\ref{thm:X_cohomology}.

With this as  background, we can state:

\setcounter{results}{1}
\begin{MainThm}\label{prop:bound_on_c1}
Let $\xi$ be a \Spinc\ structure on 4-manifold  $X$ with a fixed
point free circle action such that $SW_X(\xi)\not=0$ (in either
chamber when $b_+=1$).
\begin{enumerate}
\item If $b_+(X){>}1$ or $b_+(X){=}0$, then $c_1(\xi)$
is pulled back from $Y$.

\item If $b_+(X){=}1$, then either $c_1(\xi)$ is pulled back from $Y$, or $SW_X^+(\xi)=0$.
\end{enumerate}
\end{MainThm}

\begin{rem}
In case 2b, the Seiberg-Witten invariant of the other chamber can
be calculated using the wall crossing formula of
\cite{sw:wallcross}. \end{rem}

\begin{cor} If $b_+(X)>1$ then $c_1(\xi)^2=0$ and $X$ is SW
simple type.
\end{cor}

Recall that a 4-manifold is SW simple type if the dimension of the
moduli space is $0$ for all \Spinc\ structures with nonzero
Seiberg-Witten invariants.

Theorem A follows easily from the following formula about
Seiberg-Witten solutions.

\begin{thm}\label{thm:absolute_formula}
Let $(A,\Psi)$ be any solution in
$\mathcal{M}(X,\xi,g_X,\pi^*(\delta)^+)$.  Then
$$ 0 \ \ = \ \ \int_X |\Spnabla_{\T}\Psi|^2+ |D^+\Psi|^2
+|\iota_\T F_A|^2 + 2\pi^2 c_1(\xi)^2 + 2\pi c_1(\xi) \cdot
\pi^*(\I\delta)  .$$ The vector field $\T$ is the unit vector
field generated by the circle action.
\end{thm}

\begin{rem} The equation in Theorem \ref{thm:absolute_formula} only holds
for perturbations which are pulled back from $Y$.  It does not
hold for a general self-dual $2$-form on $X$.
\end{rem}

The rest of this subsection contains a proof of Theorem A assuming
Theorem \ref{thm:absolute_formula} above.  We will then come back
and prove Theorem \ref{thm:absolute_formula} in the next
subsection. We prove each case separately.

\bigskip

{\bf Proof of case 1:}\quad When $b_+(X)>1$, the moduli space is
nonempty for all generic metric and perturbation pairs. Since
generic pairs are dense in the space of metrics and self-dual
2-forms, we can take a sequence of generic pairs which converge
to the pair $(g_X,0)$.  By compactness, solutions of the generic
pairs converge to a solution $(A,\Psi)\in\mathcal{M}(X,\xi,g_X,0)$
and it satisfies
\begin{eqnarray}
0 \ \ = \ \ \int_X |\Spnabla_{\T}\Psi|^2+ |D^+\Psi|^2 +|\iota_\T
F_A|^2 + 2\pi^2 c_1(\xi)^2 \label{eq:terms=0}
\end{eqnarray}
by Theorem \ref{thm:absolute_formula}.  Using equation
(\ref{eq:c_1^2geq0}) we conclude that all terms in equation
(\ref{eq:terms=0}) vanish; in particular, $c_1(\xi)^2=0$ and
\begin{eqnarray*}
\iota_\T F_A &=&0.
\end{eqnarray*}
Since $dF_A=0$, this equation implies $\mathcal{L}_\T F_A=0$ by
Cartan's formula.  Together the equations $\iota_\T F_A =
\mathcal{L}_\T F_A = 0$ imply that $F_A$ is pulled back from $Y$.
Since $c_1(\xi) = \frac{\I}{2\pi}F_A$, case 1 follows.

When $b_+(X)=0$ we have that $b_2(X)=0$ is also zero because the
signature is zero.  Thus $c_1(\xi)$ is always a torsion class and
this is pulled back by Theorem \ref{thm:X_cohomology}.
\hfill$\Box$\newline

\bigskip

{\bf Proof of case 2:} \quad  Assume that $c_1(\xi)$ is not pulled
back.  By the argument proceeding the statement of Theorem A,
$c_1(\xi)^+ \not=0$.

We proceed by contradiction. Suppose that $SW_X^+(\xi)\not=0$. In
this chamber, the moduli space will be nonempty for all generic
pairs of metrics and perturbations.  Note that the unperturbed
Seiberg-Witten equations $(\delta=0)$ are in this chamber because
$(2 \pi c_1(\xi)-0)^+ \cdot c_1(\xi)^+>0$.  Hence we can use the
same argument as in case 1 to show that $c_1(\xi)$ is pulled back
from $Y$ -- contradicting our assumption.  Thus $SW_X^+(\xi)=0$.
\hfill$\Box$\newline

\subsection{Solutions to the SW equations}

In this subsection we prove Theorem \ref{thm:absolute_formula}.
The idea is to prove a Weitzenb\"{o}ck-type decomposition for the
Dirac operator we constructed in subsection \ref{sec:sw_eq_4_man}.
Before we prove this decomposition, however, we need to show that
the full Dirac operator $\dirac_A:\Gamma_X(W^+\oplus W^-) \ra
\Gamma_X(W^+\oplus W^-)$ is self-adjoint.  The following technical
lemma accomplishes this.

\begin{lem}
\label{lem:dirac_self_adjoint} Let $\xi=(W,\sigma)$ be a \Spinc\
structure over $X$.  Let  $\Spnabla$ be a spinorial connection
created by coupling a connection $A \in \mathcal{A}(\det W^+)$
with the $SO(3)$-connection $\Snabla$ defined in subsection
\ref{sec:sw_eq_4_man}. Similarly, let $\Spnabla^{\tiny L.C.}$ be
the spinorial connection created by coupling the same connection
$A$ with the Levi-Civita connection $\Snabla^{\tiny L.C.}$. Then
\begin{eqnarray}
\dirac^{\tiny L.C.}_A = \dirac_A - \frac{1}{2}\sigma(\eta \wedge
d\eta).
\end{eqnarray}
Since $\dirac^{\tiny L.C.}$ and Clifford multiplication by
$3$-forms are both self-adjoint operators, $\dirac_A$ is
self-adjoint.
\end{lem}

\begin{proof}
Extend $\eta$ to an orthonormal coframe $\{\eta=e^0, e^1,e^2,
e^3\}$ on a patch of $X$ so that $e^0 = \eta$, and $\{e^1,
e^2,e^3\}$ are horizontal lifts of an orthonormal coframe
$\{\overline{e}^1, \overline{e}^2, \overline{e}^3\}$ on $Y$. Let
$\{e_0=\T,e_1,e_2,e_3\}$ be the dual vector fields with respect
to the metric $g_X$.

The difference $1$-form $\omega = \Snabla^{\tiny L.C.} - \Snabla
\in \Omega^1(\mathfrak{so}(T^*X))$ can be thought of as an element
in $\Omega^1(\Lambda^2T^*X)$ via the vector space isomorphism
\[i:\mathfrak{so}(T^*X) \ra \Lambda^2(T^*X) \]
defined by
\[ i(a^k_j) = \frac12 \sum_{j<k}a^j_k  e^j \wedge e^k. \]
The action of $\mathfrak{so}(T^*X)$ on the bundle $W$ is modeled
on $\sigma_{\Lambda^2} \circ i$.  Thus we can  Clifford multiply
the $\Lambda^2$ component of $\omega\in\Omega^1(\Lambda^2 T^*X)$
to get
\[ \dirac_A^{\tiny L.C.} = \dirac_A + \sigma_{\Lambda^1 \ot \Lambda^2}(\omega)  \]
where $\sigma_{\Lambda^1\ot \Lambda^2}: \Lambda^1 \ot \Lambda^2
\ra \End(W)$ is a linear map defined by
\[ \sigma_{\Lambda^1 \ot \Lambda^2} (\alpha \ot \beta) = \sigma(\alpha)\sigma(\beta)\]
for a basis element $\alpha \ot \beta \in \Lambda^1\ot
\Lambda^2$. This map can be conveniently reformulated as
\[\sigma_{\Lambda^1\ot \Lambda^2}(\alpha \ot \beta)=-\sigma(\iota_{\alpha^\flat}\beta)
+\sigma(\alpha \wedge \beta),\] where $\iota_{\alpha^\flat}$ is
contraction with the vector field which is $g_X$-dual to $\alpha$.

Let $\{\zeta_{12},\zeta_{13}, \zeta_{23} \}$ be the functions
defined by
\begin{eqnarray}
\label{eq:deta}
 d \eta = 2\zeta_{12}e^1\wedge e^2 + 2\zeta_{13}e^1 \wedge e^3 + 2\zeta_{23}e^2 \wedge e^3.
\end{eqnarray}
We can use equation (\ref{eq:deta}) and the first Cartan
Structure equation
$$de^i = \sum_{j} e^j\wedge w^i_j $$
to calculate the connection matrix for $\Snabla^{\tiny L.C.}$.
For example, we can write $d\eta$ as
$$d\eta = e^1 \wedge (\zeta_{12}e^2 +\zeta_{13}e^3) +
e^2 \wedge (-\zeta_{12}e^1 +\zeta_{13}e^3) + e^3 \wedge
(-\zeta_{13} e^1 - \zeta_{23}e^2)$$ to get the top row of the
connection matrix
\begin{eqnarray}
\label{eq:LCconnectionmatrix} \left(\begin{array}{cccc} 0 &
\zeta_{12}e^2 + \zeta_{13}e^3 & -\zeta_{12}e^1 + \zeta_{23}
e^3 & -\zeta_{13}e^1 - \zeta_{23}e^2 \\
-\zeta_{12}e^2 - \zeta_{13}e^3 & 0 & -\zeta_{12}e^0 + \omega^1_2 & -\zeta_{13}e^0 + \omega^1_3\\
\zeta_{12}e^1 - \zeta_{23}e^3 & \zeta_{12}e^0 - \omega^1_2 & 0 & -\zeta_{23}e^0 + \omega^2_3 \\
\zeta_{13}e^1 + \zeta_{23}e^2 & \zeta_{13}e^0 - \omega^1_3 &
\zeta_{23}e^0 - \omega^2_3 & 0
\end{array}\right)
\end{eqnarray}
\mbox{}\\
The $\omega^i_j$'s in the second, third, and forth row are
pulled-back from the connection $1$-form for
the Levi-Civita connection on $Y$.  The connection matrix for $\Snabla$ is \\
\begin{eqnarray}
\label{eq:PBconnectionmatrix} \left(\begin{array}{cccc}
0 & 0 & 0 & 0 \\
0 & 0 &  \omega^1_2 & \omega^1_3\\
0 &  - \omega^1_2 & 0 &  \omega^2_3 \\
0 &  - \omega^1_3 &  - \omega^2_3 & 0
\end{array}\right).
\end{eqnarray}
\mbox{}\\
Using the isomorphism $i$, the difference $\Snabla^{\tiny
L.C.}-\Snabla$ can be written as
\[\omega = \frac12 \sum_{i=1}^3 e^i \ot \eta \wedge \iota_{e_i}(d\eta) + \frac12 \eta\ot d\eta. \]
A straight forward calculation gives $\sigma_{\Lambda^1\ot
\Lambda^2}(\omega) = -\frac12 \sigma(\eta\wedge d\eta)$.
\end{proof}

\begin{rem}The operators $\dirac_A^{\tiny L.C.}$ and $\dirac_A$ have the
same index.
\end{rem}

\begin{lem}
\label{lem:square} The square of the Dirac operator decomposes
into
\begin{eqnarray}
(\dirac_A^+)^* \dirac_A^+ & = & -(\Spnabla_{\T})^2 + (D^+)^* D^+
+ \frac12\sigma((\eta \wedge \iota_{\T} F_A)^+).
\end{eqnarray}
where $(\;)^+$ is the projection onto self-dual $2$-forms.
\end{lem}

\begin{proof}
We work with the full Dirac operator first.  By using the
definition of $\Snabla$ from equation (\ref{eq:_snabladef}), we
see from
\begin{eqnarray*}
\langle \sigma(\eta) \Spnabla_{\T} \Psi, \Phi\rangle
& = & \langle\Psi, \Spnabla_{\T}(\sigma(\eta)\Phi)\rangle \\
& = & \langle\Psi, \sigma(\Snabla_{\T} \eta)\Phi\rangle
+\langle\Psi, \sigma(\eta) \Spnabla_{\T}
\Phi\rangle\\
& = & \langle\Psi, \sigma(\eta) \Spnabla_{\T}\Phi\rangle
\end{eqnarray*}
that $\sigma(\eta)\Spnabla_{\T} $ is $L^2$ self-adjoint.

The Dirac operator decomposes into a sum of two self-adjoint
operators:
\[\dirac_A = \sigma(\eta) \Spnabla_{\T} + D, \]
where $D = \sum_{i=1}^3 \sigma(e^i) \Spnabla_{e_i}$.

Squaring and noting that $\Snabla_T\eta = 0$ and
$\sigma(\eta)\sigma(\eta) = -Id$ yields
\[\dirac_A^2 = -(\Spnabla_{\T})^2 + D^2 + \{\sigma(\eta) \Spnabla_{\T}, D\} \]
where $\{\cdot,\cdot\}$ denotes the anticommutator.  The last term
simplifies using Clifford relations and the equations
(\ref{eq:spinorial}) and (\ref{eq:_snabladef}):
\begin{eqnarray*}
\{\sigma(\eta) \Spnabla_{\T}, D\} &=& \sum_{i=1}^3 \sigma(\eta
\wedge e^i) [\Spnabla_{\T}, \Spnabla_{e_i}].\label{eq:liebracket}
\end{eqnarray*}
One can use the connection matrix (\ref{eq:LCconnectionmatrix})
to calculate that
\begin{eqnarray}
\label{eq:liebracket2}[\T, e_i] = \Snabla^{\tiny L.C.}_\T e_i -
\Snabla^{\tiny L.C.}_{e_i} T = 0
\end{eqnarray}
for $i=1,2,3$.   Therefore the curvature reduces to
\begin{eqnarray}
F_\Spnabla (\T,e_i) = [\Spnabla_{\T},\Spnabla_{e_i}],
\label{eq:Form_of_curvature}
\end{eqnarray}
and we can see that
\begin{eqnarray}
\{\sigma(\eta) \Spnabla_{\T}, D\} &=& \sum_{i=1}^3 \sigma(\eta
\wedge e^i) F_\Spnabla (\T, e_i) = \sigma(\eta)\sigma(\iota_\T
F_\Spnabla ). \label{eq:curvature_eq_from_dirac}
\end{eqnarray}
By the definition of $\Snabla$, the action of
$[\Spnabla_{\T},\Spnabla_{e_i}]$ for $i=1,2,3$ commutes with
Clifford multiplication.  Therefore $F_\Spnabla (\T,e_i)$ is a
scalar endomorphism, so
\begin{eqnarray}
\label{eq:curvatures_are_equal} F_\Spnabla (\T,e_i) = \frac12
F_A(\T,e_i).
\end{eqnarray}
Restricting attention to $W^+$ gives the formula.
\end{proof}

{\bf Proof of Theorem \ref{thm:absolute_formula}:}\quad Take a
solution $(A,\Psi)\in \mathcal{M}(X,g_X,\pi^*(\delta)^+)$, apply
$(\dirac_A^+)^*\dirac_A^+$, and take the inner product with
$\Psi$.  After applying Lemma~\ref{lem:square} and integrating
over $X$ we get

\begin{eqnarray}
0 &=&   \int_X \langle(\dirac_A^+)^*\dirac_A^+ \Psi, \Psi\rangle \nonumber \\
& = & \int_X \langle(\sigma(\eta)\Spnabla_{\T})^2\Psi,\Psi\rangle
+ \langle (D^+)^*D^+\Psi,\Psi\rangle + \frac12 \langle
\sigma((\eta \wedge \iota_{\T} F_A)^+)\Psi,\Psi\rangle
\nonumber \\
& = & \int_X |\Spnabla_{\T}\Psi|^2+ |D^+\Psi|^2 +  2\langle (\eta
\wedge \iota_{\T} F_A), q(\Psi) \rangle.\label{eq:main_int}
\end{eqnarray}
In the last step we used the adjoint property of $q(\Psi)$ from
equation (\ref{eq:adjointprop}) and the fact that $q(\Psi)$ is
self-dual.

Substituting the $q(\Psi)=F^+_A+\pi^*(\delta)^+$ from
(\ref{eq:sw_4dim}), we get
\begin{eqnarray}
2\langle (\eta \wedge \iota_{\T} F_A), q(\Psi) \rangle & = &
2\langle (\eta \wedge \iota_{\T} F_A),
F^+_A + \pi^*(\delta)^+) \rangle \nonumber\\
&=& \langle \iota_{\T} F_A, \;\; \iota_{\T}\displaystyle \left(
F_A +
 \star \, F_A  + \pi^*(\delta) + \star \, \pi^*(\delta)\right) \rangle \nonumber \\
&=& |\iota_{\T}F_A|^2 + \langle \iota_{\T} F_A, \iota_{\T} \star
\, \left(F_A + \pi^*(\delta)\right)
\rangle \nonumber \\
&=& |\iota_{\T}F_A|^2 + \frac12 \I F_A \wedge \I F_A  + \I F_A
\wedge \pi^*(\I\delta) \label{eq:show_piece_nonzero}
\end{eqnarray}
The last equality is true by the following calculation. Let
$F_{ij}$ be the functions defined by $$F_A = \sum_{0\leq i<j} \I
F_{ij} e^i \wedge e^j.$$ Then
$$\iota_{\T}F_A = \I F_{01}e^1 +
\I F_{02}e^2 + \I F_{03}e^3$$ and
$$\iota_{\T}\star \,F_A = \I F_{12}e^3
-\I F_{13}e^2 + \I F_{23}e^1.$$ Taking their inner product gives
$$\langle \iota_{\T}F_A,
\iota_{\T}\star \, F_A\rangle = F_{01}F_{23} -F_{02}F_{13} +
F_{03}F_{12} = \frac12 \I F_A \wedge \I F_A.$$  A similar
calculation shows that $\langle \iota_{\T}F_A, \; \iota_{\T}
\star \, \pi^*(\delta) \rangle = \I F_A \wedge \pi^*(\I\delta)$.

Integrating equation (\ref{eq:show_piece_nonzero}) over $X$ gives
the lemma.\hfill$\Box$\newline

\section{Diffeomorphic moduli spaces}
\label{sec:diff_mod}

In this section we prove Theorem B.  We continue to work with the
circle bundle $\pi:X\ra Y$ with the $S^1$-invariant metric
$g_X=\eta \ot \eta + \pi^*(g_Y)$ as in \ref{sec:sw_eq_4_man}. Fix
a closed orbifold $2$-form $\delta\in \Omega^2(Y;\I\BR)$ to
perturb the $3$-dimensional equations and pull it back to get an
$S^1$-invariant 2-form on $X$.  Perturb the 4-dimensional
equations by projecting $\pi^*(\delta)$ onto the self-dual
2-forms to get $\pi^*(\delta)^+$.

The total moduli space $\mathcal{M}(X,g_X,\pi^*(\delta)^+)$ is a
disjoint collection of moduli spaces, one component for each
\Spinc\ structure $\xi$ on $X$.  Define
$$\mathcal{N}(X,g_X,\pi^*(\delta)^+)$$ to be the components
of the total moduli space whose cohomology class $c_1(\xi)$ is
pulled back from $Y$.

Theorem A implies that we need only look at these components to
calculate the Seiberg-Witten invariants when $b_+>1$.  This
restriction on the total moduli space is done to rule out \Spinc\
structures covered in Theorem A case 2b.

Note that for any \Spinc\ structure $\xi$ whose $c_1(\xi)$ class
is pulled back and for any 2-form $\delta$,
\begin{eqnarray}
c_1(\xi)^2 &=&0 \nonumber \\
c_1(\xi)\cdot \pi^*(\I\delta) &=&0. \label{eq:N_condition}
\end{eqnarray}
In particular, the expected dimension of the moduli space is 0 by
Proposition \ref{ref:dim_=c_1^2}.

A pullback of a solution $(A_0,\Psi_0)$ to (\ref{eq:sw_3dim}) on
$Y$ is the solution $(A,\Psi) = \pi^*(A_0,\Psi_0)$ to
(\ref{eq:sw_4dim}) on $X$. Pick an orthonormal coframe on a patch
of $X$
\[\{e^0, e^1,e^2, e^3\}\]
so that $e^0 = \eta$, and $\{e^1, e^2,e^3\}$ are horizontal lifts
of an orthonormal coframe $\{\overline{e}^1, \overline{e}^2,
\overline{e}^3\}$ on $Y$. Let $\{e_0=\T,e_1,e_2,e_3\}$ be the
dual vector fields with respect to the metric $g_X$. In this case
the Dirac operator can be written as
\begin{eqnarray}
\dirac^+_{A} = \sigma(\eta)\Spnabla_{\T} + D_{A}
\label{eq:dirac_decomp}
\end{eqnarray}
where $\Spnabla$ is a connection on $W^+$ created by coupling $A$
with the reducible connection $\Snabla$ (see subsection
\ref{sec:sw_eq_4_man}) and $D=\sum \sigma(e^i) \Spnabla_{e_i}$ for
$i=1,2,3$.  From the construction of the pulled back \Spinc\
structure, it is immediately clear that $\pi^*(\Psi)$ is harmonic
since it is constant along the fiber and comes from a harmonic
spinor on $Y$.  The first equation of (\ref{eq:sw_4dim}) is
satisfied by pulling up the first equation and projecting each
term onto the self-dual $2$-forms. Since a gauge transformation
on $\xi$ pulls back to a gauge transformation of $\pi^*(\xi)$, we
get a well defined map on the level of moduli spaces.

\begin{MainThm} The pullback map $\pi^*$ induces a homeomorphism
\[ \pi^*:\mathcal{M}^*(Y,g_Y,\delta) \ra \mathcal{N}^*(X,g_X,\pi^*(\delta)^+).\]
Furthermore, if either of the two moduli spaces is a smooth
manifold, then both of them are smooth, and $\pi^*$ is a
diffeomorphism. \label{thm:diff_of_moduli_spaces}
\end{MainThm}

One remark: There is no restriction on $b_+(X)$ in the above
theorem.

\medskip

The next three subsections contains the proof of this theorem.  We
show that $\pi^*$ is a homeomorphism in the first two subsections.
In the final subsection we show that $d\pi^*$ is an isomorphism on
the kernel of the linearizations. This is sufficient to prove
that the moduli spaces are diffeomorphic because the expected
dimension of each is zero. Bochner vanishing arguments are used
to prove that $\pi^*$ and $d\pi^*$ are surjective.

\subsection{$\pi^*$ is injective}

Suppose we have two irreducible solutions to the 3-dimensional
equations whose pullbacks $(A,\Psi)$ and $(A',\Psi')$ differ by a
gauge transformation $g\in Map(X,S^1)$,
\[ g (A,\Psi) = (A',\Psi').\]
We wish to show that $g$ is in fact pulled back from $Map(Y,S^1)$.
We think of $g$ as a section of $\End(\det W^+) = \End(\pi^*(\det
W)) = \pi^*(\End(\det W))$.  Use $A$ to create a connection
$\nabla^{\End}$ on $\End(\pi^*(W^+))$ which has trivial fiberwise
holonomy. Then
\begin{eqnarray*}
(\nabla^{\End}_{\T}g)\Psi & = &
\nabla^{A}_{\T}(g\Psi) - g \nabla^{A}_{\T}\Psi \\
& = & 0
\end{eqnarray*}
because $\Psi' = g\Psi$ are pulled back sections.  By the unique
continuation theorem for elliptic operators, $\Psi \not = 0$ on a
dense open set, hence
\[\nabla^{\End}_{\T}g \equiv 0\]
on $X$.  Thus $g$ is a fiberwise constant section of the line
bundle $\pi^* (\End(\det W))$ and by Proposition~\ref{prop:1to1},
it can be pushed forward to a section of $\End(\det W)$ on $Y$,
i.e., a gauge transformation on $Y$.

\subsection{$\pi^*$ is surjective}

Take a solution $(A,\Psi)\in \mathcal{N}^*(X,g_X,\pi^*(\delta)^+)$
to the Seiberg-Witten equations (\ref{eq:sw_4dim}). We will show
that the solution is pulled up from a solution $(A_0,\Psi_0)$ on
$Y$.

Combining the formula in Theorem \ref{thm:absolute_formula} with
the fact that \Spinc\ structures from
$\mathcal{N}^*(X,g_X,\pi^*(\delta)^+)$ satisfy equations
(\ref{eq:N_condition}), we get

\begin{eqnarray}\label{eq:main_eqn_in_surjective}
0 \ \ = \ \ \int_X |\Spnabla_{\T}\Psi|^2+ |D^+\Psi|^2 +|\iota_\T
F_A|^2 + \underbrace{2\pi^2 c_1(\xi)^2}_{\mbox{\normalsize 0}}+
\underbrace{2\pi c_1(\xi) \cdot
\pi^*(\I\delta)}_{\mbox{\normalsize 0}}.
\end{eqnarray}
and that the following terms must be identically zero:
\begin{eqnarray}
\label{eq:circ_invar}
0 & = & \Spnabla_{\T}\Psi, \\
\label{eq:harmonic_3spinor} 0& = & D^+\Psi, \\
\label{eq:F_pulled_back} 0 & = & \iota_{\T}F_A,
\end{eqnarray}

Equation (\ref{eq:F_pulled_back}) implies $\mathcal{L}_\T F_A =
0$ and together the equations imply that $F_A$ is circle invariant
and pulled up from $Y$. Equation (\ref{eq:circ_invar}) and the
fact that $\Psi \not \equiv 0$ means that $A$ has trivial
fiberwise holonomy. Therefore we can apply Proposition
\ref{prop:1to1} and Remark \ref{rem:push_spinc} to $\xi$ with
connection $A$ to conclude that $\Psi$ corresponds to an orbifold
section $\Psi_0$ on a \Spinc\ structure $\xi_0$ with connection
$A_0$ on $Y$.

In this situation, $D^+$ is the Dirac operator on the orbifold
$Y$, so by equations (\ref{eq:harmonic_3spinor}) and
(\ref{eq:circ_invar}), the second Seiberg-Witten equation of
(\ref{eq:sw_3dim}) is satisfied for $(A_0,\Psi_0)$.  It is easy
to check that $(A_0,\Psi_0)$ also satisfies the first
Seiberg-Witten equation.

Therefore the map $\pi^*$ is a homeomorphism of moduli spaces.
\hfill$\Box$\newline

\bigskip

\subsection{The kernels are isomorphic}

\label{sec:ker_are_iso} Consider an irreducible solution
$\C=(A,\Psi)$ to the 4-dimensional Seiberg-Witten equations for a
fixed metric and perturbation $(g_X,\pi^*(\delta)^+)$ in the
\Spinc\ structure $\xi$. In the previous subsection we saw that
$(A,\Psi)$ was pulled back from a solution $\C_0= (A_0,\Psi_0)$
to the 3-dimensional equations on $Y$ in the \Spinc\ structure
$\xi_0$.

We now describe the tangent space at the solution $\C$.  The
following sequence of operators (for a fixed $k \geq 5$)
\[\xymatrix{
0 \ar[r] & T_{\bf 1} L^2_{k+2}(Y,S^1) \ar[r]^-{\mathcal{L}_\C} &
T_\C L^2_{k+1}(\I T^*Y \oplus W^+ ) \ar[r]^-{LSW^4} &
L^2_{k}(\I\Lambda_+ T^*Y \oplus W^- ) \ar[r] & 0}\] is called the
{\em deformation complex at $\C$}.  The map $\mathcal{L}_\C$ is
the infinitesimal action of the gauge group at $\C$ described by
its differential at the identity:
$$\mathcal{L}_\C: \I f \mapsto (-2\I df, \I f\Psi),$$
and $LSW^4$ is the linearization of the 4-dimensional
Seiberg-Witten equations with fixed perturbation
$\pi^*(\delta)^+$. We can wrap $LSW^4$ and $\mathcal{L}_\C$ into
one operator
$$\mathcal{T}_\C: L^2_{k+1}(\I T^*X \oplus W^+)\ra L^2_k(\I \Lambda_+ T^* X \oplus W^- \oplus \I \Lambda^0 T^* Y) $$ by setting
$\mathcal{T}_\C=LSW^4 + \mathcal{L}_\C^*$. Then the $\ker
\mathcal{T}_\C$ is the set of $(a,\psi)$ which satisfy
\begin{eqnarray}
\label{eq:ker_of_T_C}
\begin{array}{r@{\ = \ }@{\ 0}l}
d^+a-q(\psi,\Psi)-q(\Psi,\psi) & ,\\
\dirac_A\psi
+\frac12 \sigma(a)\Psi & ,\\
-2d^*a+\I \Imaginary \langle \psi, \Psi \rangle  &. \end{array}
\end{eqnarray}
The last equation is a slice condition for the gauge group action.

Let $\mathcal{H}_\C^*$ denote the cohomology of the deformation
complex at $\C$.  We can now state Lemma 2.2.11 from \cite[page
129]{sw:notes_on_sw_theory}:

\begin{lem}The deformation complex at $\C$ is Fredholm, that is, the
coboundary maps have closed ranges and the cohomology spaces are
finite dimensional.  Moreover,
$$\mathcal{H}^0_\C \cong \ker \mathcal{L}_\C, \;\;
\mathcal{H}^1_\C \cong \ker \mathcal{T}_\C$$ and
$$\mbox{\rm coker}\; \mathcal{T}_\C \cong \mathcal{H}^0_\C \oplus
\mathcal{H}^2_\C.$$ In particular, the expected dimension of the
moduli space for the \Spinc\ structure $\xi$ is
$$d(\xi) =
-\chi_\BR(\mathcal{H}^*_\C) = - \dim \mathcal{H}^0_\C + \dim
\mathcal{H}^1_\C - \dim \mathcal{H}^2_\C
.$$\label{lem:tangent_space}
\end{lem}

A metric and perturbation $(g,\delta)$ is called a {\em good pair}
if $\mathcal{H}^0_\C=\mathcal{H}^2_\C=0$ for every solution to
the Seiberg-Witten equations in the \Spinc\ structure $\xi$. If
$(g_X,\pi^*(\delta)^+)$ is good, then the moduli space
$\mathcal{M}(X,\xi,g_X,\pi^*(\delta)^+)$ is a smooth manifold of
dimension $d(\xi)$, its formal tangent space can be identified
with $\mathcal{H}^1_\C$ at the point $\C$, and we can use it to
calculate the Seiberg-Witten invariants of $\xi$. For a careful
treatment of these ideas, see pages 127-135 of
\cite{sw:notes_on_sw_theory}.

There is a similar complex for the solution $\C_0$ on $Y$.  It too
can be described by an operator $$\mathcal{T}_{\C_0}:
L^2_{k+1}(\I T^*Y \oplus W) \ra L^2_{k}(\I \Lambda^2 T^* Y \oplus
W\oplus \I \Lambda^0 T^* Y)$$ given by the map
\[
\left[\begin{array}{c}  a_0 \\ \psi_0\end{array}\right]
\stackrel{\mathcal{T}_{\C_0}}{\ra} \left[ \begin{array}{c}
d^+a_0-\tau(\psi,\Psi)-\tau(\Psi,\psi)\\
\dirac_{A_0}\psi_0
+\frac12 \rho(a_0)\Psi_0\\
-2d^*a_0+\I \Imaginary \langle \psi_0, \Psi_0 \rangle \end{array}
\right].\] It also has a complex at $\C_0$, and an associated
cohomology denoted by $\mathcal{H}^*_{\C_0}$ which can be
described using $\mathcal{T}_{\C_0}$ and a similar statement as
Lemma \ref{lem:tangent_space} above.

By definition $\C$ is irreducible if and only if
$\mathcal{H}^0_\C=0$ (and likewise for $\C_0$).  Hence solutions
in $\mathcal{N}^*(X,\xi,g_X,\pi^*(\delta)^+)$ satisfy
\begin{eqnarray}
0 = d(\xi) = \dim \mathcal{H}_\C^1 - \dim \mathcal{H}^2_\C
\label{eq:H1_equals_H2}
\end{eqnarray}
by Proposition \ref{ref:dim_=c_1^2},  equation
(\ref{eq:N_condition}), and the previous lemma. Therefore
$\mathcal{H}^2_\C$ vanishes for these solutions precisely when
$\dim \mathcal{H}^1_\C = 0$.  We will use this fact and the
following theorem to show when $\dim \mathcal{H}^2_\C = 0$.

\begin{thm}\label{thm:kernels_are_iso}
Let $\C=(A,\Psi)\in \mathcal{N}^*(X,\xi,g_X,\pi^*(\delta)^+)$ be
a irreducible solution to the Seiberg-Witten equations and let
$\C_0=(A_0,\Psi_0) \in \mathcal{M}^*(Y,\xi_0,g_Y,\delta)$ be the
solution such that $\C=\pi^*(\C_0)$.  Then
$$\pi^*(\mathcal{H}^1_{\C_0}) = \mathcal{H}^1_{\C}, $$
i.e., the kernels of $\mathcal{T}_{\C_0}$ and $\mathcal{T}_\C$ are
naturally isomorphic via $\pi^*$.
\end{thm}

Because the expected dimension of the moduli space on the
3-manifold is always zero, when there is a good pair
$(g_Y,\delta)$ such that $\dim \mathcal{H}^0_{\C_0} = \dim
\mathcal{H}^2_{\C_0}=\dim \mathcal{H}^1_{\C_0}=0$ for all
solutions in $\mathcal{M}(Y,g_Y,\delta)$ we get by Theorem
\ref{thm:kernels_are_iso} that the dimension of
$\mathcal{H}^1_{\C}$ will  be zero for the pulled back solutions
as well. Hence $\mathcal{H}_{\C}^2=0$ by equation
(\ref{eq:H1_equals_H2}) for all irreducible solutions $\C$
implying  that $\mathcal{N}^*(X,g_X,\pi^*(\delta)^+)$ is a smooth
manifold. Thus Theorem \ref{thm:kernels_are_iso} finishes the
proof of Theorem B. If, in addition,
$\mathcal{N}(X,g_X,\pi^*(\delta)^+)$ does not contain any
reducible solutions, then $(g_X,\pi^*(\delta)^+)$ will be a good
pair for any \Spinc\ structure pulled back from $Y$.

The rest of this subsection contains the proof of Theorem
\ref{thm:kernels_are_iso}.  We use a Bochner vanishing argument
similar to equation (\ref{eq:main_eqn_in_surjective}).

Certainly a solution to $\mathcal{T}_{\C_0}(a_0,\psi_0){=}0$ pulls
back to a solution of
$\mathcal{T}_{\C}(\pi^*(a_0),\pi^*(\psi_0))=0$. We need to show
that $\pi^*$ is surjective, i.e., for  each solution $(a,\psi)$
of the equations (\ref{eq:ker_of_T_C}), we will prove that
$$\Spnabla_{\T}\psi = 0 \qquad \mbox{and} \qquad a\in
\pi^*(\Omega^1(Y;\I\BR)).$$

Use $\eta$ to decompose $a$ into $a= f\eta + c$ where $f\in
\Omega^0(X;\I\BR)$ and $c\in \Omega^1(X;\I\BR)$. Since $(a,\psi)$
satisfies $\dirac_A^+\psi + \frac12 \sigma(a)\Psi = 0$, we have
\begin{eqnarray}
0 &=& \int_X |\dirac_A^+\psi + \frac12\sigma(a) \Psi|^2 \nonumber \\
&=& \int_X |(\sigma(\eta) \left(\Spnabla_{\T} \psi +
\frac12f\Psi)\right) + (D^+\psi + \frac12 \sigma(c) \Psi)|^2 \nonumber \\
\label{eq:ontoker} &=& \int_X |\Spnabla_{\T}\psi + \frac12 f
\Psi|^2 + |D^+\psi + \frac12 \sigma(c) \Psi|^2 + \label{eq:nonnegsum} \\
&& \hspace{2cm} 2\mbox{Re} \langle\sigma(\eta) \Spnabla_{\T}
\psi,D^+\psi\rangle + \mbox{Re}\langle\sigma(\eta) \Spnabla_{\T}
\psi, \sigma(c)
\Psi\rangle + \nonumber \\
&& \hspace{2cm} \mbox{Re} \langle f\sigma(\eta) \Psi,
D^+\psi\rangle +
 \mbox{Re} \langle  f\sigma(\eta) \Psi, \frac12
\sigma(c) \Psi\rangle. \nonumber
\end{eqnarray}

Two of the cross terms in  equation (\ref{eq:ontoker}) are zero
as follows.  First, since $\Snabla_\T\eta=0$ we have

\begin{eqnarray*}
2 \int_X \mbox{Re} \langle\sigma(\eta) \Spnabla_{\T}
\psi,D^+\psi\rangle &=& \int_X\langle\sigma(\eta) \Spnabla_{\T}
\psi,D^+\psi\rangle + \langle D^+\psi,\sigma(\eta)
\Spnabla_{\T} \psi\rangle \\
&=& \int_X \langle \psi,\sigma(\eta) \Spnabla_{\T}
(D^+\psi)\rangle + \langle\psi,D^-\left(\sigma(\eta)
\Spnabla_{\T} \psi\right)\rangle \\
&=& \int_X \langle \psi,\{\sigma(\eta) \Spnabla_{\T}
,D\}\psi\rangle
\end{eqnarray*}
But by equations (\ref{eq:curvature_eq_from_dirac}),
(\ref{eq:curvatures_are_equal}), and (\ref{eq:F_pulled_back}),
$$\{\sigma(\eta)\Spnabla_\T,D\} \ = \  \frac12 \eta \wedge \iota_\T F_A \
= \ 0.$$ Similarly, we can use the fact that $f\eta$ and $c$ are
both self-adjoint to show
\begin{eqnarray*}
2\int_X \mbox{Re} \langle \sigma(f\eta)\Psi,\sigma(c)\Psi\rangle
& =& \int_X \langle (\sigma(c)\sigma(f\eta) +
\sigma(f\eta)\sigma(c) \Psi,\Psi\rangle \\ & =&  -2 \int_X
\langle c,f\eta\rangle |\Psi|^2 \ = \ 0.
\end{eqnarray*}

The remaining two cross terms in equation (\ref{eq:ontoker}) are
analyzed in the following lemma.

\begin{lem} In the situation above,
\begin{eqnarray}
\label{eq:crossterm5} \int_X \mbox{Re}\langle\sigma(\eta)
\Spnabla_{\T} \psi, \sigma(c) \Psi\rangle +  \mbox{Re} \langle
f\sigma(\eta) \Psi, D^+\psi\rangle = \int_X |\iota_{\T}da|^2.
\end{eqnarray}
\end{lem}

\begin{proof}
Let $\{\eta=e^0,e^1,e^2,e^3\}$ be a local coframe where
$e^1,e^2,e^3$ are pulled back from the base.  First we take the
adjoints of both terms on the left hand side of equation
(\ref{eq:crossterm5}).

Applying the adjoint of $\sigma(\eta) \Spnabla_{\T}$ in the first
term of equation (\ref{eq:crossterm5}) gives
\begin{eqnarray}
\sigma(\eta) \Spnabla_{\T}(\sigma(c) \Psi) & =& \sum_{i=1}^3
\sigma(\eta) \sigma(\Snabla_{\T}(c_ie^i)) \Psi + \sigma(\eta)
 \sigma(c) (
\Spnabla_{\T}\Psi) \nonumber \\
& =& \sum_{i=1}^3 \sigma(\eta) \sigma(\Snabla_{\T}(c_ie^i)) \Psi
\nonumber \\
&=& \sum_{i=1}^3 \sigma\left((\T c_i) \eta\wedge e^i\right) \Psi
+ \sigma (\eta) \sigma(c_i\Snabla_{\T}(e^i)) \Psi
\nonumber \\
&=& \sum_{i=1}^3  \sigma\left((\T c_i) \eta\wedge e^i\right) \Psi
\label{eq:ida_part2}
\end{eqnarray}
We used equation (\ref{eq:spinorial}) in the first line, and
(\ref{eq:circ_invar}) in the second.  We also used the definition
of $\Snabla$ from equation (\ref{eq:_snabladef}).

Similarly, we take the adjoint of $D$ in the second term of
equation (\ref{eq:crossterm5}) to find
\begin{eqnarray}
\label{eq:ida_part3}
 D(f\sigma(\eta) \Psi) = \sigma(df\wedge
\eta)\Psi + \sigma(f\eta)D\Psi =\sigma(df\wedge \eta)\Psi.
\end{eqnarray}

Next we show that the sum of the right hand sides of equations
(\ref{eq:ida_part2}) and (\ref{eq:ida_part3}) is equal to
$$\sigma(\eta\wedge\iota_{\T}(da))\Psi.$$
First, note that for $i=0,1,2,3$,
\begin{eqnarray}\label{eq:i_T_of_d_eta_is_0}
\iota_T d\eta = 0 \qquad \mbox{ and } \qquad \iota_{\T} de^i = 0.
\end{eqnarray}
This holds for $e^0=\eta$ since $d\eta$ is the curvature of a
principal orbifold circle bundle so is pulled back from $Y$; it
holds for the remaining $i$ since $e^1,e^2,e^3$ are pulled back
from $Y$. Hence,
\begin{eqnarray}
 \eta\wedge
\iota_{\T} (da) &=& \eta\wedge \iota_{\T} (d (f\eta +c))
\label{eq:eta_wedge_iota_da} \\
&=& \eta\wedge \iota_{\T} \left((df \wedge \eta + f d\eta) +
\sum_{i=1}^{3}(dc_i \wedge e^i +
c_i\wedge de^i)\right) \nonumber \\
&=& df\wedge \eta + \sum_{i=1}^{3}(\T c_i )\eta\wedge e^i
\nonumber
\end{eqnarray}

Combining equations (\ref{eq:ida_part2}), (\ref{eq:ida_part3}),
and (\ref{eq:eta_wedge_iota_da}) and projecting onto the
self-dual 2-forms we get:
\begin{eqnarray*}
\int_X \mbox{Re}\langle\sigma(\eta) \Spnabla_{\T} \psi,
\sigma(c)\Psi\rangle +  \mbox{Re} \langle f\sigma(\eta) \Psi,
D^+\psi\rangle &=& \int_X \mbox{Re} \langle\psi,(\eta \wedge
\iota_{\T}da)^+ \Psi\rangle
\end{eqnarray*}
Using equation (\ref{eq:ker_of_T_C}), we can reduce further
\begin{eqnarray*}
\int_X \mbox{Re} \langle\psi, (\eta \wedge \iota_{\T}da)^+
\Psi\rangle &=& \int_X \frac12 \langle\psi, (\eta \wedge
\iota_{\T}da)^+ \Psi\rangle + \frac12 \langle(\eta \wedge
\iota_{\T}da)^+ \Psi, \psi\rangle \\
&=& \int_X 2\langle (\eta \wedge \iota_{\T}da)^+,
q(\Psi,\psi)+q(\psi,\Psi)\rangle \\
&=& \int_X 2\langle (\eta \wedge \iota_{\T}da),
d^+a\rangle \\
&=& \int_X |\iota_{\T}da|^2 + \frac12 \int_X \I da\wedge \I da.
\end{eqnarray*}
The last equality is the same calculation as in equation
(\ref{eq:show_piece_nonzero}).
\end{proof}

Combining equations (\ref{eq:nonnegsum}-\ref{eq:crossterm5}),
gives the sum of non-negative terms.  Hence we conclude that the
following terms are identically zero:
\begin{eqnarray}
\Spnabla_{\T} \psi +
\frac12f\Psi &=&0, \label{eq:temp_ker}\\
D^+\psi + \frac12 \sigma(c) \Psi &=& 0,\label{eq:3dim_ker}\\
\iota_{\T}da &=&0\label{eq:ker_i_da}.
\end{eqnarray}

Notice that equation (\ref{eq:ker_i_da}) is equivalent to
\begin{eqnarray}
\label{eq:nabla_f_equals_Tc_i} \Snabla_T a = df.
\end{eqnarray}

We investigate equation (\ref{eq:temp_ker}) more carefully in the
next lemma.

\begin{lem}
\begin{eqnarray*}
\int_X | \Spnabla_{\T} \psi + \frac12f \Psi|^2 &=& \int_X |
\Spnabla_{\T} \psi|^2 + \frac14 f^2|\Psi|^2 + 2|df|^2.
\end{eqnarray*} \label{lem:pos_ker_results}
\end{lem}

Since $f=\iota_\T a$ and $\Psi\not=0$ almost everywhere, we
conclude from this lemma that
\begin{eqnarray}
\Spnabla_{\T} \psi &=& 0,\label{eq:ker_cir_invar}\\
\iota_{\T}a&=&0.\label{eq:ker_i_a}
\end{eqnarray}

Equation (\ref{eq:ker_cir_invar}) implies that the spinor is
circle invariant while equations (\ref{eq:ker_i_da}) and
(\ref{eq:ker_i_a}) imply that $a$ is pulled back from $Y$. These
two facts together imply that $(a,\psi)$ is pulled back from some
$(a_0,\psi_0)$ on $Y$. Equation (\ref{eq:3dim_ker}) shows that
$(a_0,\psi_0)$ satisfies the last equation of
$\mathcal{T}_{\C_0}$. It is easy to verify that $(a_0,\psi_0)$
satisfies the other two equations of $\mathcal{T}_{\C_0}$. Hence
$(a,\psi)$ is in $\pi^*(\ker \mathcal{T}_{\C_0})$ and this
completes the proof of Theorem B.

\bigskip

{\bf Proof of Lemma \ref{lem:pos_ker_results}:}\quad We must show
that the cross term satisfies
\[
\int_X \mbox{Re}\langle\Spnabla_{\T} \psi, f\Psi\rangle = \int_X
2|df|^2. \label{eq:crossterm3}
\]

Integrating by parts and noting that $\Spnabla_\T\Psi = 0$,
\begin{eqnarray*}
 \int_X \mbox{Re}\langle
\Spnabla_{\T} \psi, f \Psi\rangle &=& \int_X \mbox{Re}\langle \psi, (-\Snabla_\T f) \Psi\rangle.\\
\end{eqnarray*}
Pulling out the imaginary valued function $\Snabla_\T f$, using
equation (\ref{eq:ker_of_T_C}), and integrating by parts again,
\begin{eqnarray}
\int_X \mbox{Re}\langle \psi, (-\Snabla_\T f) \Psi\rangle
\quad = \quad -2\int_X \langle \Snabla_\T f, d^*a\rangle \nonumber \\
\quad = \quad 2 \int_X \langle f, \Snabla_\T d^*a\rangle.
\label{eq:snabla_T_d_a_is_0}
\end{eqnarray}

The results follows once we show  $\displaystyle \Snabla_\T d^* a
= \Delta f$. We first calculate $d^*a$ at a point $p\in X$. Then
\begin{eqnarray*}
d^*a \quad = \quad  -\sum_{i=0}^3 \iota_{e_i} \Snabla^{L.C.}a
\quad = \quad -\sum_{i=0}^3 \langle
\Snabla_{e_i}^{L.C.}a,e^i\rangle.
\end{eqnarray*}
Differentiating this with respect to  $\Snabla_\T$,
\begin{eqnarray}
\Snabla_\T d^*a & = & -\sum_{i=0}^3 \Snabla_\T\langle
\Snabla^{L.C.}_{e_i} a, e^i \rangle \nonumber \\
&=& -\sum_{i=0}^3 \langle \Snabla_\T\Snabla^{L.C.}_{e_i} a, e^i
\rangle.\label{eq:grad_Td^*a}
\end{eqnarray}

Next we will show using the connection matrices
(\ref{eq:LCconnectionmatrix}) and (\ref{eq:PBconnectionmatrix}),
and equation (\ref{eq:liebracket2}) that
\begin{eqnarray} \label{eq:commuting_deriviatives} \sum_{i=0}^3
\langle [\Snabla_\T, \Snabla^{L.C.}_{e_i}]a,e^i \rangle &=& 0.
\end{eqnarray}
By setting $a=\sum a_k e^k$ and using the fact that $\Snabla_\T
e^i=0$,
\begin{eqnarray*}
\sum_{i=0}^3 \langle [\Snabla_\T,\Snabla^{L.C.}_{e_i}]
a,e^i\rangle & = & \sum_{i=0}^3 \langle \Snabla_\T\left((e_i
\cdot a_k) e^k + a_k \Snabla^{L.C.}_{e_i} e^k\right) -
\Snabla^{L.C.}_{e_i} \left( (\T \cdot a_k) e^k\right),
e^i\rangle \\
&=& \sum_{i=0}^3 \T\cdot e_i \cdot a_i + (\T \cdot a_k)
\langle\Snabla^{L.C.}_{e_i}e^k,e^i\rangle +  a_k \langle
\Snabla_\T \Snabla^{L.C.}_{e_i}e^k,e^i \rangle\\
&& \hspace{1.5cm} - e_i \cdot \T \cdot a_i - (\T \cdot a_k)
\langle\Snabla^{L.C.}_{e_i}e^k,e^i\rangle.
\end{eqnarray*}
The first and fourth term cancel because $[\T,e_i]=0$ by equation
(\ref{eq:liebracket2}). The second and last term also cancel. The
third term is equal to
\begin{eqnarray}\label{eq:Tonmetric}
a_k \T \cdot \langle \Snabla^{L.C.}_{e_i}e^k,e^i\rangle
\end{eqnarray}
because $\Snabla_T$ is compatible with the metric and $\Snabla_\T
e^i =0$. But
$$\sum_{i=0}^3 \langle\Snabla^{L.C.}_{e_i} e^k, e^i \rangle \quad =
\quad -\sum_{i=0}^3 \langle e^k, \Snabla^{L.C.}_{e_i}
e^i\rangle.$$  By inspecting the connection matrix
(\ref{eq:LCconnectionmatrix}), we can see that
$\Snabla^{L.C.}_{e_i}e^i$ contains only the terms $\omega^i_j$
which are pulled back from $Y$.  Since these terms are invariant
under the circle action, equation (\ref{eq:Tonmetric}) vanishes
giving equation (\ref{eq:commuting_deriviatives}).

Therefore we can commute  $\Snabla_\T$ with
$\Snabla^{L.C.}_{e_i}$ in equation (\ref{eq:grad_Td^*a}), and
apply equation (\ref{eq:nabla_f_equals_Tc_i}) to get:
\[\Snabla_\T d^*a \quad = \quad - \sum_{i=0}^3 \langle
\Snabla^{L.C.}_{e_i} df, e^i\rangle \quad = \quad  \Delta f.\]
This statement is independent of frame, so we can substitute it
into equation (\ref{eq:snabla_T_d_a_is_0}). The lemma now follows
by integration by parts.

\section{Results}

We are now ready to prove the formula for calculating the
Seiberg-Witten invariants of a $4$-manifold with a fixed point
free circle action and state some immediate corollaries.

\begin{MainThm}
\label{thm:results} Let $X$ be a closed smooth $4$-manifold with
$b_+>1$ and a fixed point free circle action.  Let $Y^3$ be the
orbifold quotient space and suppose that $\chi \in \Pic^t(Y)$ is
the orbifold Euler class  of the circle action. If $\xi$ is a
\Spinc\ structure over $X$ with $SW^4_X(\xi)\not=0$, then
$\xi=\pi^*(\xi_0)$ for some \Spinc\ structure on $Y$ and
\begin{eqnarray*}
SW^4_X(\xi) = \sum_{\xi' \equiv \xi_0 \mod \chi} SW^3_Y(\xi'),
\end{eqnarray*}
where $\xi'-\xi_0$ is a well-defined element of $\Pic^t(Y)$. When
$b_+=1$, the  formula holds for all \Spinc\ structures which are
pulled back from $Y$.
\end{MainThm}

\begin{rem}  In the $b_+(X)=1$ case, the numerical invariant may
still depend on the chamber structure of $Y$ if $b_1(Y)=1$.  A
nice example to check that Theorem A together with Theorem C
gives the correct invariants in the $b_+(X)=1$ case is the
manifold $T^2\times S^2$.
\end{rem}

\begin{proof}
Recall that for a generic choice of metric and perturbation
$(g_Y,\delta)$ the moduli space satisfies $\mathcal{H}^0_{\C_0} =
\mathcal{H}^1_{\C_0} = \mathcal{H}^2_{\C_0}=0$ for all solutions
$\C_0 = (A_0,\Psi_0)$ to the 3-dimensional Seiberg-Witten
equations (see subsection \ref{sec:ker_are_iso} for more details).
For this good pair the moduli space $\mathcal{M}(Y, g_Y, \delta)$
is a smooth manifold with no reducible solutions. Since we can
choose a perturbation generically such that the projection of
$F_{A_0} + \delta$ onto the harmonic 2-forms is not a multiple of
the harmonic representative of $\chi$ for all solutions in
$\mathcal{M}(Y,g_Y,\delta)$, we have that
\[
(\pi^*(F_{A_0}) + \pi^*(\delta))^+ \not=0
\]
on $X$ as well, hence $\mathcal{N}(X,g_X, \pi^*(\delta)^+)$ does
not contain reducible solutions either. By Theorem B,
$\mathcal{N}(X, g_X, \pi^*(\delta)^+)$ is diffeomorphic to a
smooth manifold without reducible solutions. We have in effect
shown that $(g_X,\pi^*(\delta)^+)$ is a good pair and that this
moduli space can be used to calculate the SW invariant.

Choose a specific \Spinc\ structure $\xi$ on $X$ such that
$c_1(\xi)$ is pulled back and $SW_X(\xi)\not=0$. There exists a
\Spinc\ structure $\xi_0$ on $Y$ such that $\xi=\pi^*(\xi_0)$ by
Theorem B, and
\[
\mathcal{N}(X, \xi, g_X, \pi^*(\delta)^+) = \coprod_{\xi' \cong
\xi_0 \mod_{\chi}} \mathcal{M}(Y, \xi', g_Y, \delta).
\]
From this the formula follows.\end{proof}

When the action is free, the theorem above reduces to the formula
worked out in a previous paper.

\begin{cor}
Let $X$ be a closed smooth $4$-manifold with $b_+>1$ and a free
circle action.  Then the orbit space $Y^3$ is a smooth
$3$-manifold and suppose that $\chi \in H^2(Y;\BZ)$ is the first
Chern class of the circle action on $X$.  If $\xi$ is a \Spinc\
structure over $X$ with $SW^4_X(\xi)\not=0$, then
$\xi=\pi^*(\xi_0)$ for some \Spinc\ structure on $Y$ and
\begin{eqnarray*}
SW^4_X(\xi) = \sum_{\xi' \equiv \xi_0 \mod \chi} SW^3_Y(\xi'),
\end{eqnarray*}
where $\xi'-\xi_0$ is a well-defined element of
$H^2(Y;\BZ)$.\label{sw:baldridge1}
\end{cor}

Because of this formula, it is particularly easy to calculate the
Seiberg-Witten invariants for manifolds with free circle
actions.  The proof of this theorem in \cite{sw:circleactions} is
completely different than the proof of Theorem C.  It used a
gluing theorem and the next corollary, which now is just a
consequence of Theorem C:

\begin{cor}[c.f. Donaldson \cite{sw:swand4man}]
Let $X \cong Y^3\times S^1$ with $b_+(X)>1$. If a \Spinc\
structure $\xi$ has  $SW_X(\xi)\not=0$, then there is one \Spinc\
structure $\xi_0$ on $Y$ such that $\xi=\pi^*(\xi_0)$  and
\[SW^4_X(\xi) = SW^3_Y(\xi_0).\]\label{sw:donaldson}
\end{cor}

The usual route used to explain the corollary above is to consider
the cyclic covering of $X$ by $Y^3 \times \BR$. There is a natural
way to pullback solutions of (\ref{eq:sw_4dim}) to solutions on
$Y^3 \times \BR$ for \Spinc\ structures pulled up from $Y^3$.
After putting the solution in temporal gauge it satisfies the
3-dimensional Seiberg-Witten equations because it is a constant
gradient-flow of the Chern-Simons-Dirac functional
\cite{symp:symp_lef_fib_S1xM}.  Thus for each $\xi$ on $X$ such
that $SW_X(\xi)\not=0$ there is a \Spinc\ structure on $Y$ whose
moduli space is nonempty for all generic metrics and
perturbations.  This corollary shows that this moduli space can
actually be identified with the moduli space of $X$ and can be
used to calculate the Seiberg-Witten invariant.

\section{Examples}

\label{sec:example_b+=1} In this subsection we construct  a
$b_+(X){=}1$ 4-manifold with free circle actions whose
Seiberg-Witten invariants are still diffeomorphism invariants. In
this situation we can use Theorem C to calculate its
Seiberg-Witten polynomial. We then use Theorem B to study the
moduli spaces of $X$ and its quotient $Y$ and explain why the
invariants do not change when crossing a ``wall.''

\begin{figure}[h]
\begin{center}
\includegraphics{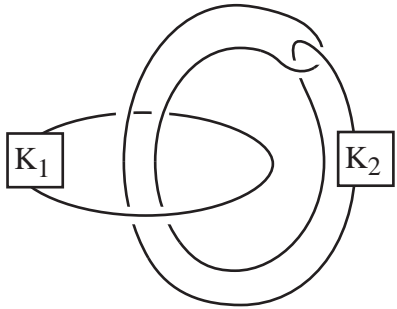}
    \caption{$Y_{K_1K_2}$ before surgery. }
    \label{fig:whitehead}
\end{center}
\end{figure}

Recall the construction from \cite{sw:circleactions}. Take the
Whitehead link in $S^3$ and compose each component with the knots
$K_1$ and $K_2$ (see Figure \ref{fig:whitehead}). Then the
3-manifold $Y_{K_1K_2}$ is the result of surgery on this new link
with each surgery coefficient equal to 0. Because the Whitehead
link is fibered, when and $K_1$ and $K_2$ are fibered knots, the
resulting 3-manifold fibers over the circle.

Define the $X_{K_1K_2}(L)$ to be the unit circle bundle of a line
bundle $L$ over $Y_{K_1K_2}$.  When $c_1(L)$ is nontorsion, we
get the following facts:

\begin{enumerate}
\item $b_1(Y_{K_1K_2})=2$, $b_1(X_{K_1K_2}(L)) = 2$, and
$b_+(X_{K_1K_2}(L))=1$.
\item The cup product pairing
$$\cup:H^1(X_{K_1K_2}(L);\BZ) \ot H^1(X_{K_1K_2}(L);\BZ)\ra H^2
(X_{K_1K_2}(L);\BZ)$$ is trivial.  This can be computed from the
cup product on $Y_{K_1K_2}$ using the isomorphism
$\pi^*:H^1(Y_{K_1K_2};\BZ) \ra H^1(X_{K_1K_2}(L);\BZ)$.
\end{enumerate}

The two facts above are exactly the conditions needed to show that
the wall crossing number is zero for all \Spinc\ structures
\cite{sw:wallcross}. Hence Seiberg-Witten invariants are still
diffeomorphism invariants for these manifolds.  In fact, any unit
circle bundle over a three manifold which satisfies the
conditions above will be such an example.  The manifolds
constructed above are also particularly easy for calculating  the
Seiberg-Witten polynomial using Theorem C.  We give one example.

\begin{figure}[h]
\begin{center}
\includegraphics{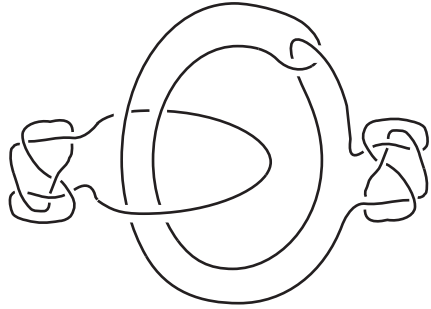}
    \caption{$Y$ constructed out of $6_3$ knots. }
    \label{fig:6_3_knot}
\end{center}
\end{figure}

Let $Y=Y_{K_1K_2}$ be the manifold where  $K_1$ and $K_2$ are the
fibered $6_3$ knot in \cite{knots:knotsandlinks} (see Figure
\ref{fig:6_3_knot}). Then the Seiberg-Witten polynomial
$$\mathcal{SW}^3_{Y}(x,y) = (x^{-4} -3x^{-2} +5 -3x^2 +
x^4)(y^{-4} -3y^{-2} +5 -3y^2 + y^4)$$ is calculated using Milnor
torsion and Meng-Taubes' theorem \cite{sw:swequalmilnortorsion}.
In this setup $x=\mbox{exp}(PD(m_1))$ and $y=\mbox{exp}(PD(m_2))$
are formal variables where $m_1, m_2 \in H_1(Y;\BZ)$ represent
the meridian loops of each component of the Whitehead link.  Thus
the term $9x^{2}y^{2}$ in the polynomial above means that the
Seiberg-Witten invariant for the \Spinc\ structure identified
with $PD(2m_1+2m_2)$ is $9$.

Let $X=X_{K_1K_2}(L)$ be the unit circle bundle of a line bundle
$L$ which satisfies $c_1(L)=4PD(m_1)$.  By Theorem A and the fact
that Seiberg-Witten invariants for $X$ are independent of the
wall crossing, $c_1(\xi)$ is pulled back from $Y$. Thus
$d(\xi)=0$ and condition (\ref{eq:N_condition}) holds for \Spinc\
structures with nontrivial SW invariants. We can apply Theorem C
to get
$$\mathcal{SW}^4_{X}(x,y) =
7y^{-4}-6x^2y^{-4}-21y^{-2}+18x^2y^{-2}+35-30x^2-21y^2+18x^2y^2+7y^4-6x^2y^4$$
where the formal variables are defined by
$x=\mbox{exp}(\pi^*(PD(m_1)))$ and $y=\mbox{exp}(\pi^*(PD(m_2)))$
and represent the pullback of \Spinc\ structures on $Y$.  Note
that $X$ is an example of a nonsymplectic 4-manifold with a circle
action whose quotient fibers over the circle.

Theorem C gives insight into why the Seiberg-Witten invariant does
not change when crossing a wall. Let $G_X$ be the product space of
metrics and $\Gamma(\Lambda_+)$ and recall that  $(g_X,\delta) \in
G_X$ is called a {\em good pair} if the moduli space
$\mathcal{M}(X,\xi,g_X,\delta)$ is a smooth manifold without
reducible solutions.  When $b_+>1$ the wall of bad pairs is at
least codimension 2 and a cobordism can be constructed between
the two moduli spaces of good pairs. However, when $b_+(X)=1$ it
is possible that two good pairs cannot be connected through a
generic smooth path in $G_X$ without crossing a wall of bad pairs
where reducible solutions occur. Passing through a bad pair could
cause a singularity to occur in the cobordism. For a general
$b_+=1$ manifold, this will often break the invariance of the
Seiberg-Witten invariant.

Suppose that we had two good pairs that can not be connected
without going through a bad pair. Connect the two good pairs with
smooth generic paths to good pairs of the form $(g_X=\eta^2 + g_Y,
\pi^*(-F \pm \eta)^+)$. Here $g_X$ is fixed, $||\eta||$ is
sufficiently small, and $F$ is the harmonic curvature form which
represents $2\pi \I c_1(\xi_0)$ for some \Spinc\ structure on $Y$.
Suppose for the sake of argument that $\{\gamma(t) = (g_X,
\pi^*(-F + t\eta)^+)\; | \; -1 \leq t\leq 1\}$ is a smooth generic
path in $G_X$ connecting the good pairs. Then a bad pair occurs
in both $G_X$ and $G_Y$ precisely when $t=0$. While the wall has
codimension $b_+(X)=1$ in $G_X$ and hence unavoidable, the wall
in $G_Y$ has codimension $b_1(Y)=2$. Thus it is possible to
perturb the path in $G_Y$ to a smooth generic path which avoids
any bad pairs. The moduli spaces $\mathcal{M}(Y,\xi_0,g_Y, -F \pm
\eta)$ are then cobordant, and
$$SW_Y(\xi_0,g_Y,-F-\eta) = SW_Y(\xi_0,g_Y,-F+\eta).$$
This can be done for each \Spinc\ structure $\xi'$ on $Y$ such
that $\xi = \pi^*(\xi')$, so by Theorem~C
$$SW_X(\xi,g_X,\pi^*(-F-\eta)^+) = SW_X(\xi,g_X,\pi^*(-F+\eta)^+), $$
i.e., the Seiberg-Witten invariant is independent of metric and
perturbation.

Note that the perturbed path in $G_Y$ will correspond to a
perturbed path in $G_X$ which will still go through a bad pair.
The moduli space for $X$ will have reducible solutions at the bad
pair, but they do not change the value of the Seiberg-Witten
invariant.

The same analysis holds for any $b_+{=}1$ 4-manifold with a fixed
point free circle action and $b_1(Y)=2$.  Therefore we get the
following corollary to Theorem C.

\begin{cor}
Let $X$ be a $b_+{=}1$ 4-manifold with a fixed point free circle
action whose quotient $Y$ satisfies $b_1(Y)=2$.  If
$\xi=\pi^*(\xi_0)$ is a \Spinc\ structure which is pulled back
from a \Spinc\ structure $\xi_0$ on $Y$, then $$SW^4_X(\xi) =
\sum_{\xi' \equiv \xi_0 \mod \chi} SW^3_Y(\xi')$$ and the
numerical invariant does not depend on the chamber in which it was
calculated.
\end{cor}



\section{Final remarks}
\label{sec:final_remarks}

Theorem C together with an affirmative answer to the following
conjecture would establish a way to calculate Seiberg-Witten
invariants for {\em any} $b_+{>}1$ 4-manifold with a circle
action.

\begin{conj}
If  $X$ is a $b_+>1$ smooth closed $4$-manifold with a circle
action that has fixed points, then $\mathcal{SW}_X \equiv 0$.
\end{conj}

There is already considerable evidence which suggest that this is
true.  For simply connected 4-manifolds carrying a circle action,
we can apply the classification result of R. Fintushel
\cite{circact:Circ_act_on_four_man, circact:Class_of_circ_4man}.

\begin{thm}[Fintushel] Modulo the 3-dimensional Poincar\'{e} conjecture, a
simply connected 4-manifold carrying a smooth $S^1$-action must
be a connected sum of copies of $S^4$, $\CP^2$,
$\overline{\CP}^2$, and $S^2\times S^2$.
\end{thm}

This classification result is enough to show that in the $b_+>1$
case, $X$ is the connected sum of two $b_+>0$ pieces, and hence
$\mathcal{SW}_X\equiv 0$.

The conjecture also follows from Proposition 4 of
\cite{sw:circleactions} for 4-manifolds with smooth semi-free
actions whose orbit space $Y$ has a nonempty boundary and the
rank $H_1(Y,\partial Y;\BZ)>1$.

A counter example to the conjecture above would be just as
interesting as the proof.






\begin{thebibliography}{99999}

\bibitem[B]{sw:circleactions} S. Baldridge, {\em Seiberg-Witten invariants of
$4$-manifolds with free circle actions}, Communications in
Contemporary Mathematics, accepted for publication, to appear.
E-print GT-9911051.

\bibitem[CM]{symp:symp_lef_fib_S1xM} W. Chen and R. Matveyev, {\em Symplectic
Lefschetz fibrations on $S^1\times M^3$}, Geom. Topol. {\bf 4}
(2000), 517--535.

\bibitem[D]{sw:swand4man} S. Donaldson, {\em The Seiberg-Witten equations
and 4-manifold topology}, Bull. A.M.S. {\bf 33} (1996), 45--70.



\bibitem[F1]{circact:Circ_act_on_four_man} R. Fintushel, {\em Circle actions
on simply connected 4-manifolds},
Trans. A.M.S. {\bf 230} (1977), 147--171.

\bibitem[F2]{circact:Class_of_circ_4man} R. Fintushel, {\em Classification of
circle actions on 4-manifolds}, Trans. A.M.S. {\bf 242} (1978),
377--390.

\bibitem[FS2]{sw:knots_links_and_four_man} R. Fintushel and R. Stern, {\em
Knots, links, and 4-manifolds}, Invent. Math., {\bf 134} (1998),
363--400.

\bibitem[FM1]{sw:alg_surface} R. Friedman and J. Morgan, {\em Algebraic
surfaces and Seiberg-Witten invariants}, J. Algebraic Geom., {\bf
6} (1997), 445--479.

\bibitem[FM2]{sw:obs_semireg_swinv} R. Friedman and J. Morgan, {\em Obstruction bundles,
semiregularity, and Seiberg-Witten invariants}, Comm. Anal. and
Geom., {\bf 7} (1999), 451--496.

\bibitem[FuS]{circact:sf_hom_3_sphere} M. Furuta and B. Steer, {\em Seifert fibred
homology 3-spheres and the Yang-Mills equations on Riemann
surfaces with marked points}, Adv. in Math., {\bf 96} (1992),
38--102.

\bibitem[HP]{circact:circ_act_on_4man_II} W. Huck and V. Puppe, {\em Circle actions
on 4-manifolds II}, Arch. Math. {\bf 71} (1998), 493--500.


\bibitem[LL]{sw:wallcross}T.J. Li and A. Liu, {\em General wall crossing formula},
Math.l Res. Letters  {\bf 2} (1995), 797-810.

\bibitem[M]{sw:sw_eq_and_fourmanifolds} J. Morgan, `The
Seiberg-Witten Equations and Applications to the Topology of
Smooth Four manifolds', Princeton University Press, Princeton,
New Jersey, 1996.

\bibitem[MT]{sw:swequalmilnortorsion} G. Meng and C. Taubes, {\em
${\underline{SW}}={\text{Milnor Torsion}}$}, Math. Research
Letters {\bf 3} (1996), 661--674.


\bibitem[MOY]{sw:sw_inv_seifert_space} T. Mrowka, P. Ozsv\'{a}th, and B. Yu, {\em
Seiberg-Witten monopoles on Seifert fibered spaces}, Comm. Anal.
Geom. {\bf 5} (1997), 685 -- 791.

\bibitem[N]{sw:notes_on_sw_theory} L. Nicolaescu, 'Notes on Seiberg-Witten Theory',
Graduate Studies in Mathematics {\bf 28}, American Mathematical
Society, Providence, Rhode Island, 2000.

\bibitem[OS1]{sw:higher_type_adj_ineq} P. Ozsv\'{a}th and Z. Szab\'{o}, {\em
Higher type adjunction inequalities in Seiberg-Witten theory},
preprint.

\bibitem[OS2]{sw:symp_thom_conj} P. Ozsv\'{a}th and Z. Szab\'{o}, {\em
The symplectic Thom conjecture}, Ann. of Math. {\bf 151} (2000),
93 -- 124.

\bibitem[R]{knots:knotsandlinks} D. Rolfsen, 'Knots and Links', Publish or
Perish, Inc., Houston, TX, 1989.

\bibitem[S]{orb:gauss_bonnet} I. Satake, {\em The Gauss-Bonnet theorem for V-manifolds},
J. Math. Soc. Japan   {\bf 4} (1957), 464--491.

\bibitem[T]{sw:sw_inv_and_symp_form} C. Taubes, {\em The Seiberg-Witten
invariants and symplectic forms}, Math. Res. Letters   {\bf 1}
(1994), 809--822.

\bibitem[W]{sw:monopolesandfourman} E. Witten, {\em Monopoles and
four-manifolds}, Math. Res.
Letters {\bf 1} (1994), 769--796.

\end{thebibliography}
\end{document}